\documentclass{article}

\usepackage{amsmath, amsthm, amssymb,graphics}

\usepackage{authblk, float, bm}

\usepackage{hyperref}
\usepackage{relsize}

\usepackage[caption = false]{subfig}
\usepackage[final]{graphicx}

\title{%
  Weighted Lagrange Interpolation Using Orthogonal Polynomials \\
  \large Stenger{'}s Conjecture, Numerical Approach}

\author[1]{Maha Youssef\thanks{maha.youssef@uni-greifswald.de: Corresponding author}}
\author[2]{Gerd Baumann\thanks{gerd.baumann@uni-ulm.de}}
\affil[1]{\small{Institute of Mathematics and Computer Science, University of Greifswald, Walther-Rathenau-Stra{\ss}e 47, 17489 Greifswald, Germany}}
\affil[2]{\small{Mathematics Department, German University in Cairo, Egypt. University of Ulm, Albert-Einstein-Allee 11, Ulm, Germany}}

\date{}

\newtheorem{thm}{Theorem}
\newtheorem*{prf}{Proof}

\newtheorem{dfn}{Definition}
\newtheorem{ex}{Example}[section]
\newtheorem{expt}{Experiment}

\usepackage{fancyhdr}
\begin{document}

\maketitle
\begin{abstract}
In this paper we investigate polynomial interpolation using orthogonal polynomials. We use weight functions associated to orthogonal
polynomials to define a weighted form of Lagrange interpolation. We introduce an upper bound of error estimation for such kind of approximations.
Later, we introduce the sufficient condition of Stenger{'}s conjecture for orthogonal polynomials and numerical verification for such conjecture.
\end{abstract}


\section{Introduction}

The solution of partial differential equations (PDEs) is based on the function approximation using polynomials, rational functions or, Sinc functions.
The target is to get accurate approximations of convolution integrals, Laplace and Fourier transform and, their inverses based on polynomials. Polynomial
approximation using orthogonal polynomials has been introduced before in \cite{Stenger_86}. In \cite{Stenger_86}, the Chebyshev polynomials have been used in conjunction with
Lagrange interpolation to give an accurate function approximation. In this paper, we introduce a weighted polynomial approximation with a weight
function \(\xi (x)\) that uses orthogonal polynomials with weight functions \(w(x)\). This weighted polynomial approximation is defined via Lagrange
interpolation using the roots of the orthogonal polynomials as interpolation point. We will show that this weighted Lagrange approximation has an
exceptional rate of convergence different from a standard polynomial approximation. The convergence rate of the integral operators (convolution,
Inverse Laplace and Inverse Fourier) is based on properties of the formed matrices. This property has been formulated by Stenger in \cite{Stenger_2015} as {``}New
Polynomial Conjecture{''}. This conjecture is proved for a restrict case, but not proven beyond this. In this paper, we re-formulate the conjecture
by an extra sufficient condition and verify this numerically for high degree orthogonal polynomials. In addition, we verify it for a new set of polynomials
defined by Lagrange interpolation at Sinc points, Poly-Sinc polynomials \cite{Stenger_2013}. On the other hand we show that the conjecture formulated in \cite{Stenger_2015} is not
always true by introducing some contradicting examples.

This paper is organized as follows, in Section 2, we introduce weighted Lagrange approximation using orthogonal polynomials and their rate of convergence.
In Section 3, we introduce Stenger{'}s conjecture and the main theorem of convergence for convolution, Laplace and, Fourier operators. In addition
we introduce our conditioned version of Stenger{'}s conjecture. In Section 4, we introduce simulation results for both conjectures. Finally, comments
and conclusion are given in section 5.


\section{Polynomial Approximation}

Consider an interpolation procedure on an interval \((a,b)\), of the form 

\begin{equation}
f(x)\approx \sum _{k=1}^n b_k(x)f\left(x_k\right),
 \label{equation:1}
\end{equation}

where \(a<x_1<x_2<\cdot \cdot \cdot <x_{n-1}<x_n<b\), and \((a,b)\subseteq \mathbb{R}\). The \(b_k\){'}s are the basis functions, which can be polynomials
or Sinc functions \cite{Stenger_2011}. For polynomial interpolation, we have

\begin{equation}
b_k(x)=\prod _{l=0, l\neq k}^n \frac{x-x_l}{x_k-x_l}=\frac{v(x)}{\left(x-x_k\right)v'\left(x_k\right)}, v(x)=\prod _{l=0}^n x-x_l.
 \label{equation:2}
\end{equation}

Lagrange interpolation is used for computing over finite intervals, \(-\infty <a<x<b<\infty\). Recently it was shown that it can be used effectively
on \((0,\infty )\) or $\mathbb{R}$ \cite{Maha_Phd}. The idea of this extension is based on the choice of the interpolation points \(x_k\). 

An alternative to Lagrange interpolation is a Sinc interpolation, we have thus in (\ref{equation:1})

\begin{equation}
b_k(x)=\frac{\sin \left\{\frac{\pi }{h}(\phi (x)-k h)\right\}}{\frac{\pi }{h}(\phi (x)-k h)},
 \label{equation:3}
\end{equation}

\noindent where \(\phi :(a,b)\longrightarrow \mathbb{R}\) is a one-to-one transformation (conformal map). As in Lagrange interpolation, Sinc can be defined
in any interval, finite or infinite, based on the proper choice of the conformal maps \(\phi (x)\). For example, in case of \((a,b)=\mathbb{R}\),
we choose \(\phi (x)=x\). Using this conformal map a set of interpolation points is created on an interval \((a,b)\) as 

\begin{equation}
x_k=\phi ^{-1}(k h),
 \label{equation:4}
\end{equation}

\noindent where \(k=-M, \text{... }, N\) and \(h\) is a step length. These points give us the flexibility to define Lagrange interpolation not only on a finite
interval, as usual, but also on \(\mathbb{R}^+=(0,\infty )\) and \(\mathbb{R}=(-\infty ,\infty )\). Of course a proper conformal map must be chosen
in each case \cite{Stenger_2011}.


\subsection{Orthogonal basis}

In the definition of the basis function \(b_k(x)\) in (\ref{equation:2}), we use the function \(v(x)=p_n(x)\) where \(p_n\){'}s, \(n=1,2,3,\ldots\) are a set of orthogonal
polynomials with weight function \(w(x)\). For instant, for the finite case, \(p_n(x)\) can be chosen as Legendre polynomials, Chebyshev polynomials,
Jacobi polynomials,..etc. For the semi-infinite case the \(p_n\){'}s are the Laguerre polynomials while in the infinite case Hermit polynomials are
used. \(p_n(x)\) have the same definition as \(v(x)\) defined in (\ref{equation:2}), with \(x_k\), the interpolation points, are the roots of the orthogonal polynomials
in the corresponding interval, i.e. solution of the equation \(p_n(x)=0\).

\begin{equation}
\begin{split}
b_k(x) & =\prod _{l=1, l\neq k}^n \frac{x-x_l}{x_k-x_l}=\frac{p_n(x)}{\left(x-x_k\right)p_n'\left(x_k\right)},\\
 p_n(x)& =\prod _{l=1}^n x-x_l\text{    }\text{ with }\text{     }p_n\left(x_k\right)=0.
\end{split}
 \label{equation:5}
\end{equation}

This approximation provides an accurate approximation for the function as well as for the integral of the function. 

Using the basis functions (\ref{equation:5}), we shall define a row vector \(\pmb{B}\) of basis functions and an operator vector \(V\) that maps a function \(f\)
into a column vector of order \(n\) by

\begin{equation}
\pmb{B}(x)=\left(b_1,\cdot \cdot \cdot ,b_n\right), \text{   }
V(f)=\left(f\left(x_1\right),\text{... }, f\left(x_n\right)\right)^T.
 \label{equation:6}
\end{equation}

This notation enables us to write the above interpolation schemes in an operator form, as

\begin{equation}
g\approx \pmb{B}\text{   }V(f).
 \label{equation:7}
\end{equation}


\subsection{Error Estimation}

\begin{dfn} 
Let D be a simply connected domain having a boundary \(\partial D\), and let \(a\) and \(b\) denote two distinct points of \(\partial D\). Let D(y,r)=\(\{z\in \mathbb{C}, \left|z-y\right|<r\}\). Define
\begin{equation*}
D_2=D\cup D(y,r), \text{for } y\in (a,b).
\end{equation*}
\end{dfn}

\begin{thm}
Let \(f\) be an analytic and bounded function in \(D_2\) and let \(r>0\). Define the Lagrange approximation of \(f\) as (\ref{equation:1}) and
(\ref{equation:5}) then there exist two constants \(A\) and \(B\), independent of \(n\) and \(r\), such that

\begin{equation}
\left\| E_n(f,x)\right\| =\left\| f-\pmb{B}\,V(f)\right\| \leq \min \left\{A r^{-n},B \sqrt{n}r^{-n}\right\},
 \label{equation:8}
\end{equation}

\end{thm}

\begin{prf}
We prove this theorem in the case of Legendre polynomials. In this case \(a=-1\) and \(b=1\). The error of \(f\) via the Lagrange approximation using
Legendre polynomials \(p_n\) can be expressed as a Cauchy Taylor contour integral

\begin{equation}
E_n(f,x)=\frac{p_n(x)}{2\pi  i}\int _{\partial D}\frac{f(z)}{(z-x)p_n(z)}dz,\text{    }z\in \partial D_2 \text{and } x\in [-1,1].
 \label{equation:9}
\end{equation}

From the definition of \(D_2\), \(\left| z-x\right| \geq r\) for all \(x\in [-1,1]\) and \(z\in \partial D_2\) and so \(p_n(z)=\prod _{l=0}^n z-x_l\geq
r^n\). As \(f\) is bounded in \(D_2\) then \(\left| f\right| \leq \mathfrak{B}(f)\) in \(D_2\). Then

\begin{equation}
\left\| E_n(f,x)\right\| =\left| f(x)-\pmb{B}\, V(f)\right| \leqslant \frac{\mathfrak{B}(f)}{r^{n+1}}\frac{L\left(\partial D_2\right)}{2\pi}\underset{x\in [-1,1]}{\max }\left\| p_n(x)\right\| ,
 \label{equation:10}
\end{equation}

where \(L\left(\partial D_2\right)\) is the length of the \(\partial D_2\) and \(L\left(\partial D_2\right)\leq c+2\pi  r\), with \(c>0\) is a positive
integer independent of \(r\). Then

\begin{equation}
\left\| E_n(f,x)\right\| =\left\| f(x)-\pmb{B}\, V(f)\right\| \leqslant A\frac{1}{r^n}\underset{x\in [-1,1]}{\max }\left\| p_n(x)\right\|.
 \label{equation:11}
\end{equation}

For the \(\underset{x\in [-1,1]}{\max }\left\| p_n(x)\right\|\), we use the Bernstein inequality \cite{Bogaert_2012, Gautschi_2009}

\begin{equation}
\left\| p_n(x)\right\| \leq \frac{2}{\sqrt{\pi (2n+1)\sqrt{1-x^2}}}.
 \label{equation:12}
\end{equation}

This sharp upper bound of Legendre polynomials is the perfect except that it has singularities at \(x=\pm 1\), for this we use the integral definition of Legendre polynomials as

\begin{equation}
\left\| p_n(x)\right\| \leq \frac{1}{\pi }\int_0^{\pi }\left(x^2+\left(1-x^2\right)\cos^2\theta\right)^{n/2} \, d\theta \leq 1.
 \label{equation:13}
\end{equation}

The integrand on the right hand side of the inequality has its maximum value, \(1\), at \(\cos ^2\theta =1\), i.e. at \(x=\pm 1\). Using (\ref{equation:12}) and (\ref{equation:13}) in (\ref{equation:11}) yields

\begin{equation*}
\left\| E_n(f,x)\right\| \leqslant \min \left\{A\frac{1}{r^n}, B \sqrt{n}r^{-n}\right\}
\end{equation*}
\noindent with \(A\) and \(B\) are positive and independent of \(n\) and \(r\).
\end{prf}

Remarks:

\begin{itemize}
\item The convergence bound \(r^{-n}\) in (\ref{equation:8}) doesn{'}t occur as \(n\to \infty\) unless \(r\geq 1\).

\item For different types of orthogonal polynomials, suitable Bernstein type inequalities should be used instead of (\ref{equation:12}) and (\ref{equation:13}) in the above proof. For
some of these Bernstein type inequalities for Jacobi, Gegenbauer, and Laguerre polynomials, see \cite{Meinardus_67, Gautschi_2009, Zhao_2013, Koornwinder_2018}.

\item The upper bound in (\ref{equation:8}) is a pessimistic estimation of the error rate. The numerical calculations show better convergence rates. 
\end{itemize}


\subsection{Weighted Lagrange approximation}

In this section we introduce a weighted version of Lagrange approximation defined in (\ref{equation:1}) and (\ref{equation:5}). This formula is given by

\begin{equation}
f(x)\approx \xi (x)\sum _{k=1}^n b_k(x)\frac{f\left(x_k\right)}{\xi \left(x_k\right)},
 \label{equation:14}
\end{equation}

where \(b_k(x)\) are the basis functions defined in (\ref{equation:5}) and \(\xi (x)\) is an arbitrary weight function.

\begin{thm}
Let \(f/\xi\) be an analytic and bounded function in \(D_2\) and let \(r>0\). Define the Lagrange approximation of \(f\) as
(\ref{equation:14}) and (\ref{equation:5}) then there exist a constant C, independent of \(n\), such that

\begin{equation}
\left\| f-\xi (x)B\, V\left(\frac{f}{\xi }\right)\right\| \leq C \left\| E_n\right\| ,
 \label{equation:15}
\end{equation}

where \(\left\| E_n\right\|\) is the error defined in (\ref{equation:8}).
\end{thm}

\begin{prf}
The proof follows from Theorem 1 and the fact that each weight function is bounded.
\end{prf}

			
\begin{ex}{\bf{Finite Case}}
\newline

Let \(f(x)=(1+x)^{\frac{1}{3}}(1-x)^{\frac{1}{2}}\) which is an analytic bounded function on \([-1,1]\). In Fig. 1, the errors convergence using
Legendre, Chebyshev, Jacobi, and Gegenbauer roots are given. For the error estimation, we used the \(L_2-\text{norm }\) error. The figure represents
the logarithm of the error fitted as linear decaying. The logarithmic plot in Fig. 1 shows that the error \(E_n\) is qualitatively following the
decaying error given in (\ref{equation:15}). In the calculations, we used \(\xi (x)=w(x)\), where \(w(x)\) is the associated weight function for each class of orthogonal
polynomials. 
\end{ex}

\begin{figure}[H]
\centering
\includegraphics[scale=1.0]{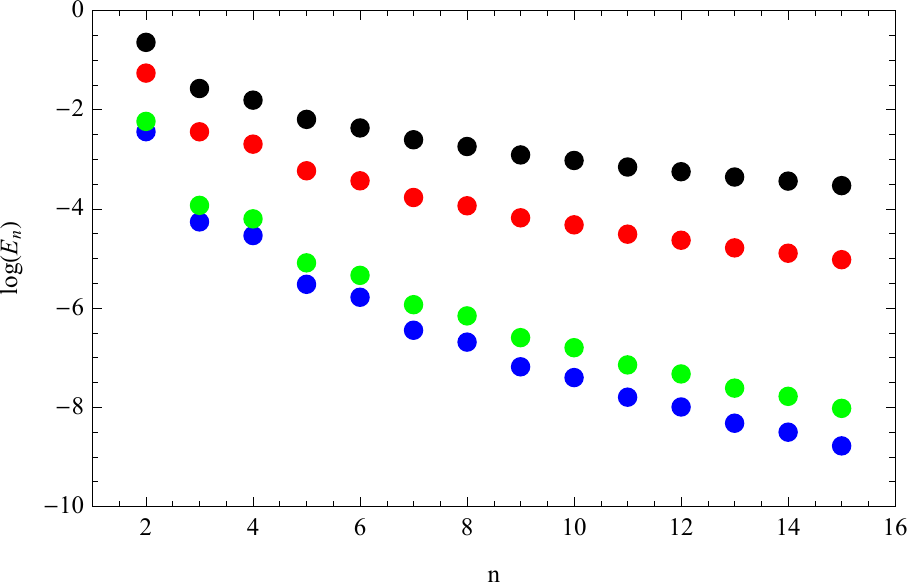}
\caption{Convergence rate of example 1 using orthogonal polynomials. Black for Chebyshev first kind polynomials \(T_n\), red for Legendre polynomials \(P_n\), green for Gegenbauer polynomials \(C_n^{(2)}\), and blue dots for Jacobi polynomials \(J_n^{(2,2)}\)}
\label{Fig:1}
\end{figure}


\begin{ex}{\bf{Semi-Infinite case}}
\newline

Let \(f(x)=x^{\frac{3}{4}}(2+x)^{\frac{1}{2}}e^{-2x}\) which is an analytic bounded function on \([0,\infty )\). In Fig. 2, the errors convergence
is shown using Laguerre polynomials and weight function \(\xi (x)=w(x)=e^{-x}\). The figure represents the logarithmic \(L_2-\text{norm }\) error.
The logarithmic plot in Fig. 2 shows that the error \(E_n\) is qualitatively following the decaying error given in (\ref{equation:15}).
\end{ex}

\begin{figure}[H]
\centering
\includegraphics[scale=1.0]{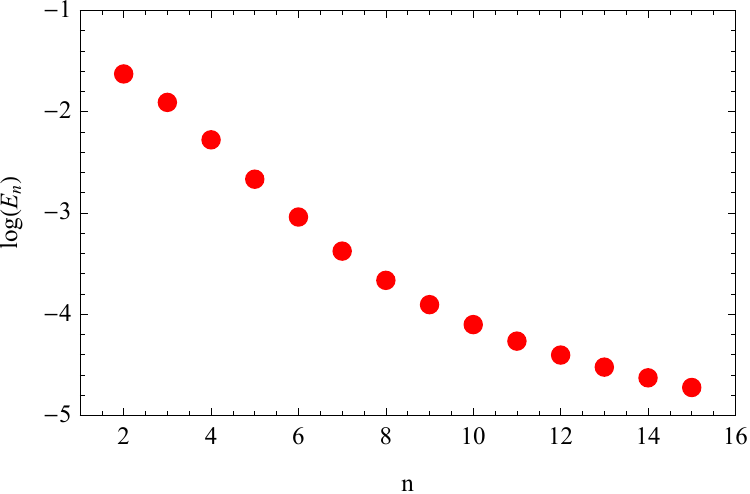}
\caption{Convergence rate of example 2 using Laguerre polynomials.}
\label{Fig:2}
\end{figure}

\newpage

\begin{ex}{\bf{Infinite case}}
\newline
Let \(f(x)=\left(x^2+2\right)^{-1/3}\left((x+2)^2+4\right)^{-1/4}\) which is an analytic bounded function on \(\mathbb{R}\). In Fig. 3, the error
convergence using Hermite polynomials and the weight function \(\xi (x)=w(x)=e^{-x^2}\) is shown. 
\end{ex}

\begin{figure}[H]
\centering
\includegraphics[scale=1.0]{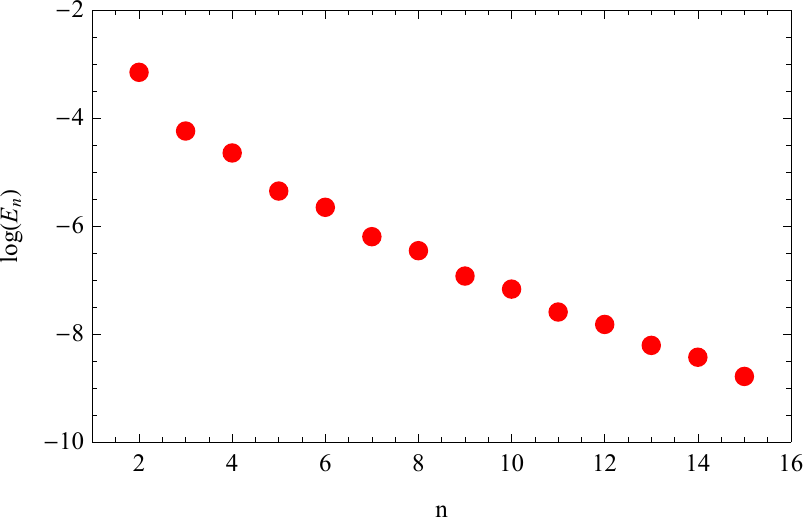}
\caption{Convergence rate of example 3 using Hermit polynomials.}
\label{Fig:3}
\end{figure}


\section{Stenger{'}s Conjecture}

In his article in 2015, Stenger conjectured \cite{Stenger_2015} that the eigenvalues of a discrete indefinite integral based on a basis \(b_k(x)\) is always positive
if \(b_k(x)\) is a polynomial. This conjecture was slightly altered by Gautschi in his proof for Legendre polynomials only \cite{Gautschi_2018}. However, we will
verify or falsify numerically the Stenger conjecture up to a large number of discretizing points which are typically not used in approximations;
i.e. we are limiting the application range of the method up to a maximal number of discretizing points.


\subsection{Indefinite Integral}

Now, define numbers \(\beta ^{\pm }\) and \(m\times m\) matrices \(B^{\pm }\) by

\begin{equation}
\begin{split}
\beta_{j k}^{+}  =\int^{x_j}_a {b_k (x)\, dx}, & \,\,\, \beta_{j k}^{-}=\int^b_{x_j} {b_k (x)\,dx} \text{  with}\\
B^{\pm  } & =\left[\beta _{j k}^{\pm }\right].
\end{split}
 \label{equation:16}
\end{equation}

The eigenvalues of the matrices \(B^{\pm }\) obey the following proportional relation: if $\lambda $ is an eigenvalue of \(B^+\) (or \(B^-\)) defined
on a finite interval \((a,b)\) then \(\tilde{\lambda }=\left.\left(\tilde{b}-\tilde{a}\right)\right/(b-a) \lambda\) is the eigenvalue of the matrix
\(\tilde{B^+}\) (or \(\tilde{B^-}\)) defined on an interval \(\left(\tilde{a},\tilde{b}\right)\). In addition, the eigenvectors of these matrices
are independent of the length of the interval. The matrices \(B^{\pm }\) and \(\tilde{B^{\pm }}\) obey the same proportionality relations. In addition,
similar simple transformations apply for infinite intervals. This enables efficient storage of a small number of such matrices. These matrices are
useful in approximating the following integrals,

\begin{equation}
\left(\mathcal{J}^+g\right)(x)=\int_a^x g(t) \, dt,\\
\left(\mathcal{J}^-g\right)(x)=\int _x^bg(t)\, dt.
 \label{equation:17}
\end{equation}

These two operators \(\mathcal{J}^{\pm }g\) are approximated by \(\mathcal{J}_m^{\pm }g\) defined as

\begin{equation}
\begin{split}
\left(\mathcal{J}^+g\right)(x) & \approx \left(\mathcal{J}_m^+g\right)(x)=\pmb{B}(x)\text{   }B^+\text{   }V (g),\\
\left(\mathcal{J}^-g\right)(x) & \approx \left(\mathcal{J}_m^-g\right)(x)=\pmb{B}(x)\text{   }B^-\text{   }V(g).
\end{split}
 \label{equation:18}
\end{equation}

We thus get the approximation \cite{Stenger_2011}

\begin{equation}
\mathcal{J}^{\pm }g\approx \mathcal{J}_m^{\pm }g.
 \label{equation:19}
\end{equation}

The matrices \(B^{\pm }\) can be explicitly expressed using Sinc quadrature that will be discussed in the next section. The integrals \(\mathcal{J}^{\pm
}\) are used to approximate both Fourier and Laplace inverse operators. These definitions play a crucial role in the solution of PDEs using indefinite
convolution representation. More precisely, the inverse Fourier transform is

\begin{equation}
f^{\mp }=\left(1\left/\mathcal{J}^{\pm }\right.\right)\tilde{f^{\mp }}\left(\mp i\left/\mathcal{J}^{\pm }\right.\right)I,
 \label{equation:20}
\end{equation}

where \(I(x)=1\) is a constant function defined on \((0,\infty )\) and { }

\begin{equation}
\tilde{f^{\mp }}(y)=\int _0^{\infty }f^{\mp }(x)e^{\mp i x y}dx,
 \label{equation:21}
\end{equation}

is the Fourier transform of the function \(f^{\mp }\in L^2(0,\infty )\). Similarly inverse Laplace transform is given by

\begin{equation}
f=\left(1\left/\mathcal{J}^+\right.\right)F^{+ }\left(\mathcal{J}^+\right) I,
 \label{equation:22}
\end{equation}

where \(F^{\mp }\) is the Laplace transform defined as

\begin{equation}
F^{\pm }(s)=\int _0^cf(\pm t)e^{-t/s}dt,
 \label{equation:23}
\end{equation}

where \((0,c)\subseteq (0,\infty )\). { }Consider the two indefinite convolution integrals on an interval \((a,b)\subseteq \mathbb{R}\),

\begin{equation}
q_a(x)=\int _a^xf(x-t)g(t)dt,\text{                           }q_b(x)=\int _x^bf(t-x)g(t)dt.
 \label{equation:24}
\end{equation}

If \(q_a\) and \(q_b\) are the convolution integral in (\ref{equation:24}), then

\begin{equation}
q_a=F\left(\mathcal{J}^+\right) V(g)\text{               }\text{and }\text{          }q_b=F\left(\mathcal{J}^-\right) V(g)
 \label{equation:25}
\end{equation}

\begin{thm}
If the spectrum of { }\(B^{\pm }\) is on the open right half plane, then there exist four constants \(C_1\), \(C_2\), \(C_3\) and,
\(C_4\) independent of \(m\) such that 

\begin{equation*}
\begin{split}
\left\| f^{\mp }-\left(1\left/\mathcal{J}_m^{\pm }\right.\right)\tilde{f^{\mp }}\left(\mp i\left/\mathcal{J}_m^{\pm }\right.\right)I\right\| & \leq C_1 \epsilon _m,\\
\left\| f-\left(1\left/\mathcal{J}_m^+\right.\right)F^{+ }\left(\mathcal{J}_m^+\right) I\right\| & \leq C_2 \epsilon _m,\\
\left\| q_a -F\left(\mathcal{J}_m^+\right) g\right\| \leq C_3 \epsilon_m\text{    }\text{and } & \left\| q_b-F\left(\mathcal{J}_m^-\right) g \right\| \leq C_4 \epsilon_m,
\end{split}
\end{equation*}

\noindent where the error \(\epsilon _m\) is a linear combination of the two errors \(\left\| \mathcal{J}^+-\mathcal{J}_m^+\right\|\) and \(\left\|f-\tilde{f}\right\|\).

\end{thm}

\begin{prf}
For the proof, see \cite{Stenger_2011}.
\end{prf}

Theorem 3 shows that if the real part of the eigenvalues of the matrices \(B^{\pm }\) are always positive then both inverse Laplace and inverse Fourier
transforms are convergent. More over, it has been shown that the rate of this convergence dependent on both function approximation and quadrature
technique used to calculate \(\mathcal{J}_m^{\pm }\). In addition, if the real parts of the eigenvalues of \(B^{\pm }\) are on the right half plane,
then the inverses of these matrices exist and can be used to yield an accurate approximation of the function derivative as

\begin{equation*}
V (f')\approx \pm \left(B^{\pm }\right)^{-1}V (f)
\end{equation*}

This eigenvalue property of \(B^{\pm }\) is defined as a conjecture by Stenger first for Sinc approximation and then for polynomial approximation.
For Sinc approximation, the conjecture has been proved recently in \cite{Han_2014}. In \cite{Stenger_2015}, Stenger introduced a conjecture related to the matrices \(B^{\pm
}\) defined by polynomial basis, \(b_k(x)\). Recently, Gautschi introduced a proof of the polynomial conjecture in a restricted special case.


\subsection{Polynomial Conjecture}

In this section we describe the new polynomial conjecture (NPC) formulated by Stenger in \cite{Stenger_2015}. 

\vspace{10pt}
{\textit{\textbf{New Polynomial Conjecture (NPC):}}} Let \(\xi (x)\) denote a function that is positive a.e. on an interval \((a,b)\subseteq \mathbb{R}\), and
assume that the moments \(M_k=\int _a^b\xi (x)x^kdx\) exists and are finite for all NonNegative integers k. Let \(\left\{p_m\right\}_{m=0}^{\infty
}\) denote the sequence of orthogonal polynomials, with \(p_m(x)\) of degree m in x, i.e.,

\begin{equation}
\int _a^b  w(x)\, p_m(x)\,p_n(x)dx=c \delta _{m,n},
 \label{equation:26}
\end{equation}

where \(c>0\) and where \(\delta _{m,n}\) denotes the Kronecker delta. If for m$\geq $1, \(p_m(x)=0\) for \(x=x_{-M}<\text{... }<x_N\), with \(m=M+N+1\),
and if the numbers \(\beta _{j k}^+\) are defined by

\begin{equation}
\begin{split}
\beta _{j k}^+ & =  \int_a^{x_j}\xi (x)\frac{p_m(x)}{\left(x-x_k\right)p_m^{'}\left(x_k\right)}dx,\\
 \beta _{j k}^- & = \int_{x_j}^b\xi (x)\frac{p_m(x)}{\left(x-x_k\right)p_m^{'}\left(x_k\right)}dx,
\end{split}
 \label{equation:27}
\end{equation}

\noindent then every eigenvalue of \(B^{\pm }=\left[\beta _{j k}^{\pm }\right]\) lies on the open right half of the complex plane.

The function \(w(x)\) in (\ref{equation:26}) is the weight function associated to the class of orthogonal polynomials \(p_n(x)\). While the function \(\xi (x)\)
is a chosen function such that the moments \(M_k\) are finite. Generally, \(w(x)\neq \xi (x)\). In Gautschi \cite{Gautschi_2018}, this conjecture has been proved
for Legendre with \(w(x)=\xi (x)=1\) and for a special class of Jacobi polynomials with \(w(x)=1-x\) and using \(\xi (x)=1\). In this paper, we first
discuss the sufficient condition for the choice of the function \(\xi (x)\). Then, we introduce a numerical verification of the new conditioned conjecture,
showing that with the sufficient condition \(\xi (x)=w(x)\), the conjecture is always verified. We test orthogonal polynomials defined on finite,
semi-infinite and infinite intervals. In addition, we test for orthogonal set of polynomials defined at Sinc points. Later, we discuss the case of
\(\xi (x)=1\) for all the orthogonal polynomials. The case of \(\xi (x)=1\), will show that the conjecture is not verified for all classes of orthogonal
polynomials.

\vspace{10pt}
{\textit{\textbf{Conditioned Polynomial Conjecture (CPC):}}} Let \(\xi (x)\) denote a function that is positive a.e. on an interval \((a,b)\subseteq \mathbb{R}\),
and assume that the moments \(\int _a^b\xi (x)x^kdx\) exists and are finite for all NonNegative integers k. Let \(\left\{p_m\right\}_{m=0}^{\infty
}\) denote the sequence of orthogonal polynomials, with \(p_m(x)\) of degree m in x, i.e.,

\begin{equation*}
\int _a^b w(x)\, p_m(x)\, p_n(x)\, dx=c \, \delta _{m,n},
\end{equation*}

\noindent where \(c>0\) and where \(\delta _{m,n}\) denotes the Kronecker delta. If for n$\geq $1, \(p_n(x)=0\) for \(x=x_0<\text{... }<x_n\) and if \(\xi(x)=w(x)\) then every eigenvalue of \(B^{\pm }=\left[\beta _{j k}^{\pm }\right]\) lies on the open right half of the complex plane, with \(\beta_{j k}^{\pm }\) defined in (\ref{equation:27}).

\vspace{8pt}
The conditioned polynomial conjecture (CPC) gives a sufficient condition for the NPC to be true. If the functions \(\xi (x)\) in (\ref{equation:27}) is the weight
function, corresponding to each class of orthogonal polynomial, then the polynomial conjecture is always true.


\subsection{Sinc Quadrature}

To compute the matrices \(B^{\pm } = \left[\beta _{j k}^{\pm }\right]\), we need to calculate the integrals defined in (\ref{equation:27}). One of the most efficient
techniques is the one based on Sinc methods \cite{Stenger_2011}. 

Given the interpolation in (\ref{equation:1}) with basis functions defined in (\ref{equation:3}), the integration

\begin{equation}
\int_a^b f(x) \, dx\approx h V(1/\phi ') V(f),
 \label{equation:28}
\end{equation}

where \(\phi (x)\) is a conformal map from \([a,b]\) onto $\mathbb{R}$ and where both vectors \(V(1/\phi ')\) and \(V(f)\) are calculated at Sinc
points defined as

\begin{equation}
x_k=\phi ^{-1}(k h), \text{with } k=-N,\text{... },N \text{and } h=\frac{\pi }{\sqrt{N}}.
 \label{equation:29}
\end{equation}

 The formula (\ref{equation:28}) is simply the Trapezoidal rule after applying the conformal map and that \(x\) is replaced by \(\phi ^{-1}\) \cite{Stenger_2011} . The approximation
in (\ref{equation:28}) has an error that is exponentially decaying with the number of Sinc basis used in the approximation,

\begin{thm} \textsl{Sinc quadrature} \cite{Stenger_2011}

\begin{equation}
\left|\int_a^b f(x) \, dx-h V(1/\phi ') V(f)\right|\leq C \sqrt{N} \text{Exp }\left(-\sqrt{\alpha  \pi  N}\right),
 \label{equation:30}
\end{equation}

where $\alpha $ is a positive constant and C is a constant independent of N.
\end{thm}


\section{Simulation Results}

In this section, we test and verify the CPC for different families of orthogonal polynomials. Some of these polynomials are defined on finite intervals
while the others are defined over semi-infinite or infinite intervals. In addition, we verify the conjecture for the basis of Poly-Sinc approximation,
which is Lagrange approximation in connection with conformal maps, that is covering the three cases of intervals. For the finite interval case, we
use the following theorem,

\begin{thm} 
If the function \(\xi (x)\) is symmetric on \([-a,a]\) then the spectrum of { }\(B^+\) is the same as the spectrum of \(B^-\).
\end{thm}

\begin{prf}
See \cite{Gautschi_2018}.
\end{prf}

\subsection{Verification of CPC}

In this section, we test numerically the CPC. First, we verify it for a set of orthogonal polynomials defined on finite interval. Specifically, we
consider, Legendre polynomials, Chebyshev polynomials of first and second kind, Jacobi polynomials and Gegenbauer polynomials. Second, we verify
it for orthogonal polynomials defined on semi-infinite and infinite intervals. Specifically, Laguerre and Hermit polynomials.

\begin{expt}{\textbf{CPC for Legendre Polynomials}}

Legendre polynomials \(p_n(x), n=0,1,2,\text{... }\) and \(x\in [-1,1]\) are orthogonal via the weight function \(\xi (x)=w(x)=1\). We can verify
the following conditions of the conjecture above (see Fig. 4), that

\begin{equation*}
M_k=\int_{-1}^1 1\times  x^k \, dx\geq 0, \forall k\geq 1.
\end{equation*}

\end{expt}

\begin{figure}[H]
\centering
\includegraphics[scale=1.0]{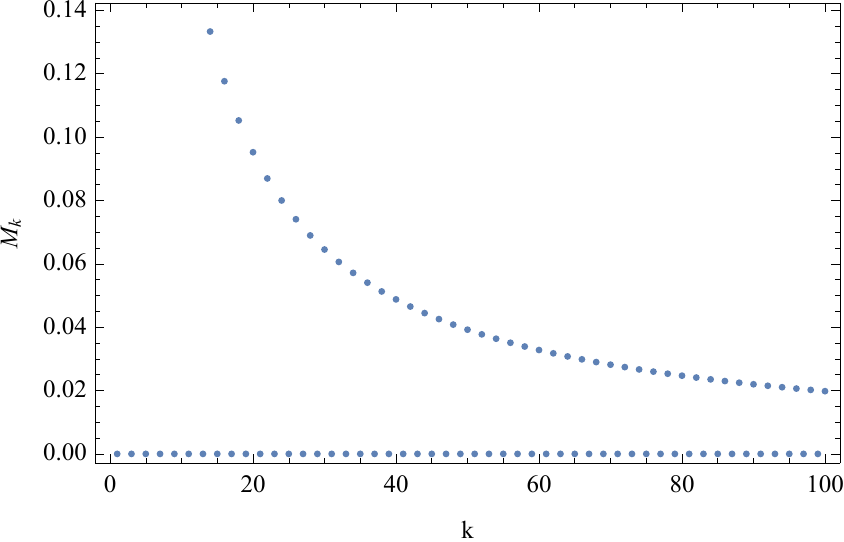}
\caption{The condition \(M_k\) for \(k=\)1,2,...,100 using \(\xi (x)=1\). The zero values corresponding to odd values of \(k\).}
\label{Fig:4}
\end{figure}

In this case we have \(\xi (x)=w(x)=1\). We are interested of two matrices \(B^{\pm }=\left[\beta _{j k}^{\pm }\right]\), with \(\beta _{j k}^{\pm}\) are defined in (\ref{equation:27}). The calculations in (\ref{equation:27}) are done using Sinc quadrature defined above with a conformal map defined on finite intervals \(\left[-1,x_j\right]\) and \(\left[x_j,1\right]\). Finally we run \(j\) from \(1\) to \(n\) to get the \(n\times n\) matrices \(B^{\pm }.\) Fig. 5 represents the eigenvalues of the matrix \(B^+\) using Legendre polynomials of degree \(n=100, 200, \dots, 500\). Note, higher degree \(n\) results to the smaller values of the absolute value of both real and imaginary part of the eigenvalues.

\begin{figure}[H]
\centering
\includegraphics[scale=1.0]{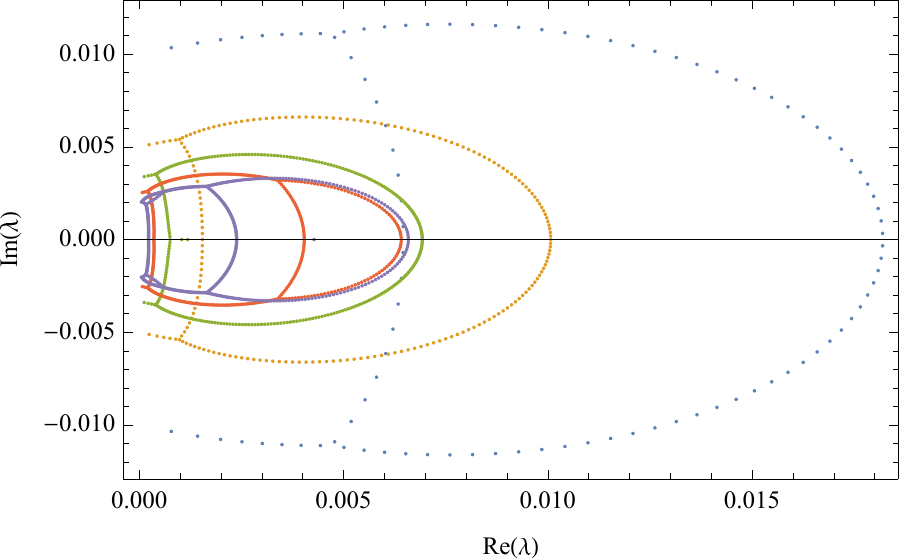}
\caption{The eigenvalues of \(B^+\)for Legendre polynomials using \(n=100, 200, 300, 400, 500\).}
\label{Fig:5}
\end{figure}

To discuss the spectral properties of the matrix \(B^+\) (or \(B^-\)) we can use the properties of the resolvent of \(B^+\) (or \(B^-\)). For the norm of the resolvent, \(\left\| \left(B^+ -\lambda  I\right)^{-1}\right\| =1\left/\text{dist }\left(\lambda ,\sigma \left(B^+\right)\right)\right.\), where \(\sigma \left(B^+\right)\) is the spectrum of \(B^+\), we can detect the distribution and magnitude of the eigenvalues. Some of the calculations
for different \(n\) are given in Fig. 6. 

\begin{figure}
\subfloat{\includegraphics[width=0.45\linewidth]{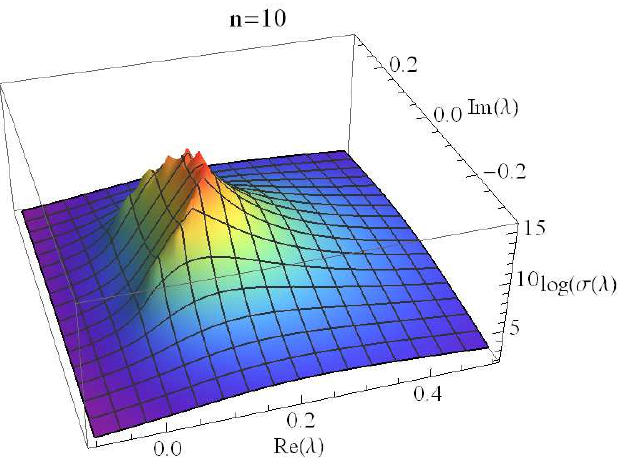}} 
\subfloat{\includegraphics[width=0.4\linewidth]{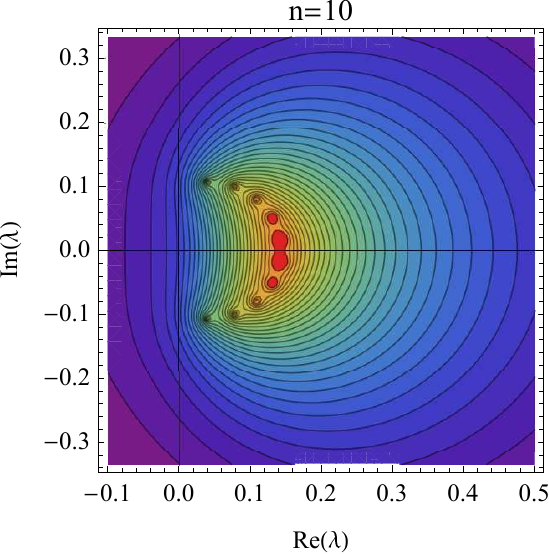}}\\
\subfloat{\includegraphics[width=0.45\linewidth]{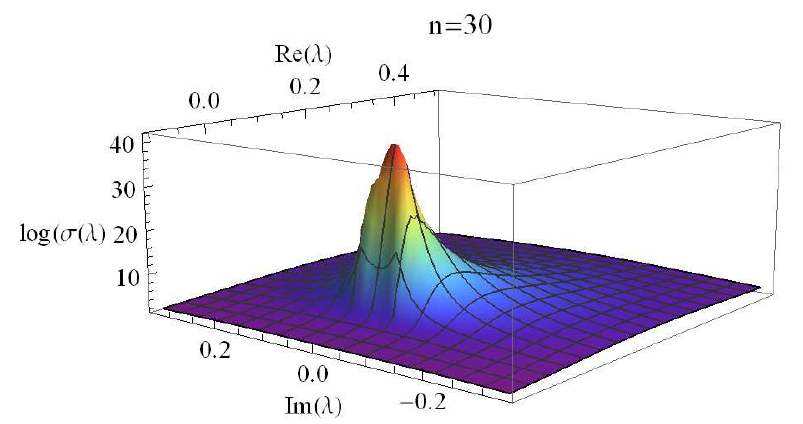}}
\subfloat{\includegraphics[width=0.4\linewidth]{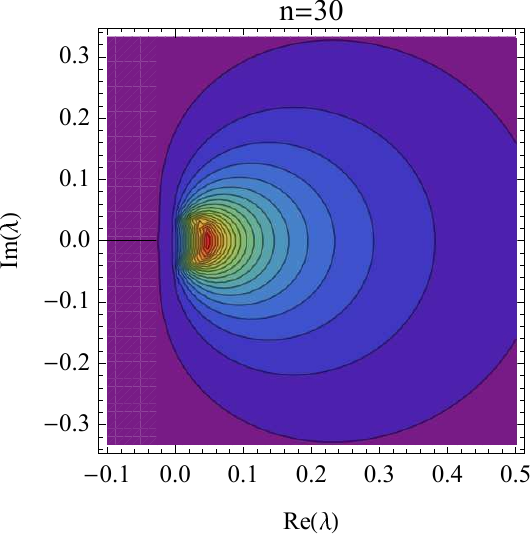}}\\
\subfloat{\includegraphics[width=0.45\linewidth]{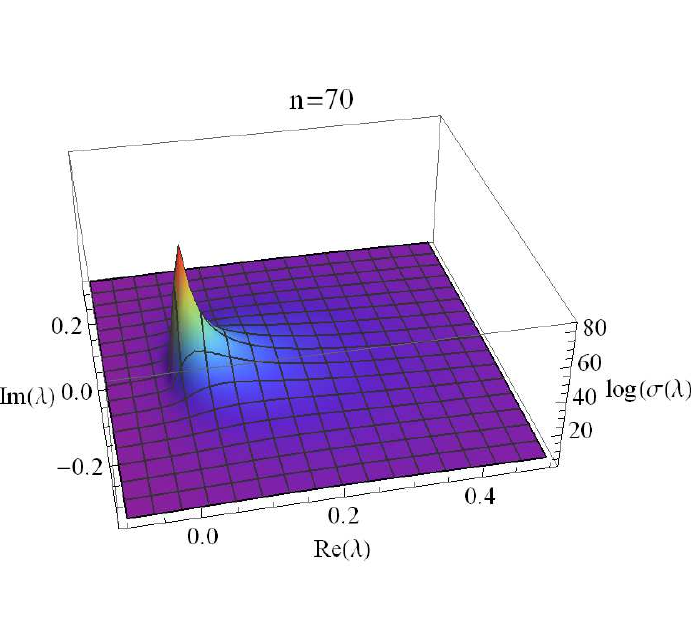}}
\subfloat{\includegraphics[width=0.4\linewidth]{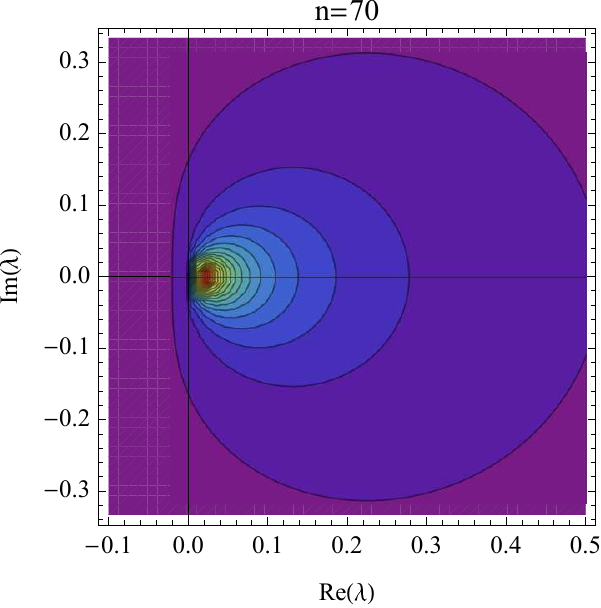}}
\caption{Resolvent of \(B^+\) for Legendre polynomials.}
\label{Fig:6-11}
\end{figure}


\begin{expt}{\textbf{ CPC for Chebyshev Polynomials}}

We define Chebyshev polynomials of the first kind, \(T_n(x)\), and second kind, \(U_n(x)\), with weight functions \(\xi (x)=w(x)=\left(1-x^2\right)^{\mp
1/2 }\), respectively. We can verify the conjecture conditions (see Figure 7), that

\begin{equation*}
M_k=\int _{-1}^1x^k \left(1-x^2\right)^{\mp 1/2} dx\geq 0, \forall k\geq 1.
\end{equation*}

\begin{figure}[H]
\centering
\includegraphics[scale=1.0]{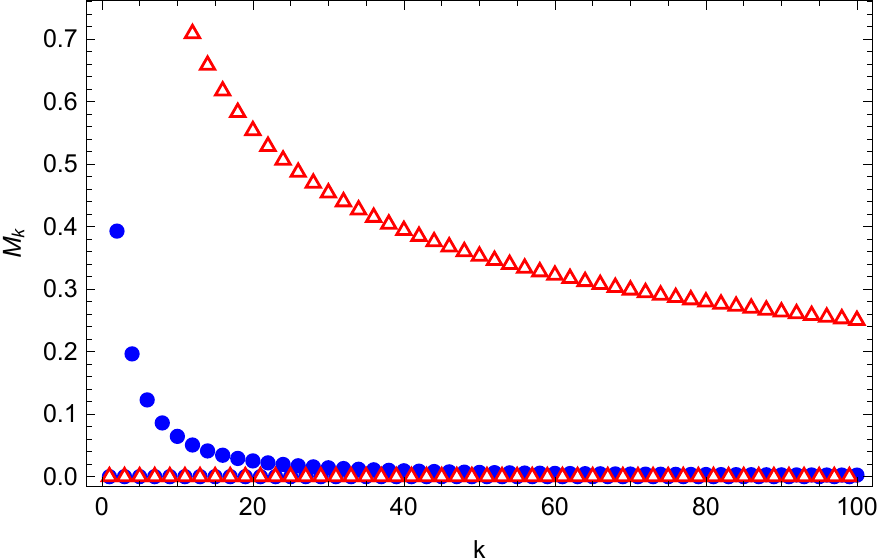}
\caption{The condition \(M_k\) for \(k=\)1,2,...,100 using \(\xi (x)=\left(1-x^2\right)^{\mp \frac{1}{2}}\).}
\label{Fig:12}
\end{figure}

As \(\xi (x)=\left(1-x^2\right)^{\mp 1/2}\) is symmetric, then we introduce here the eigenvalue computations of \(B^+\) only, these calculations are given in Fig. 8. The  resolvent is shown in Fig. 9 for different \(n\). 
\end{expt}

\begin{figure}[H]
\subfloat{\includegraphics[width=0.55\linewidth]{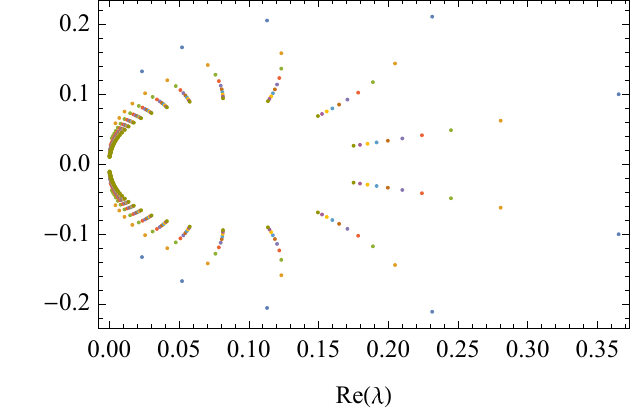}} 
\subfloat{\includegraphics[width=0.55\linewidth]{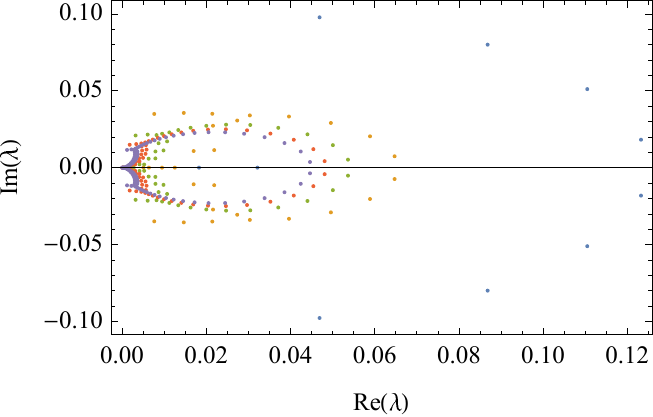}}
\caption{The eigenvalues of \(B^+\)for \(T_n(x)\), left panel, and for \(U_n(x)\), right panel. \(n=10,20, \ldots,100\).}
\label{Fig:13-14}
\end{figure}

\begin{figure}
\subfloat{\includegraphics[width=0.45\linewidth]{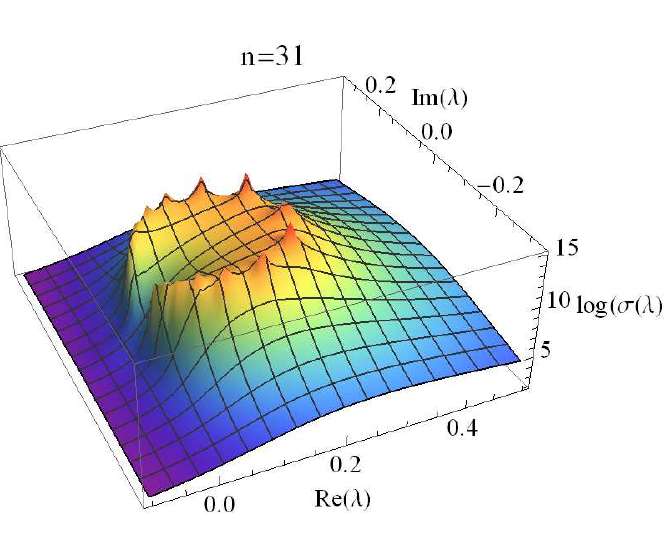}} 
\subfloat{\includegraphics[width=0.4\linewidth]{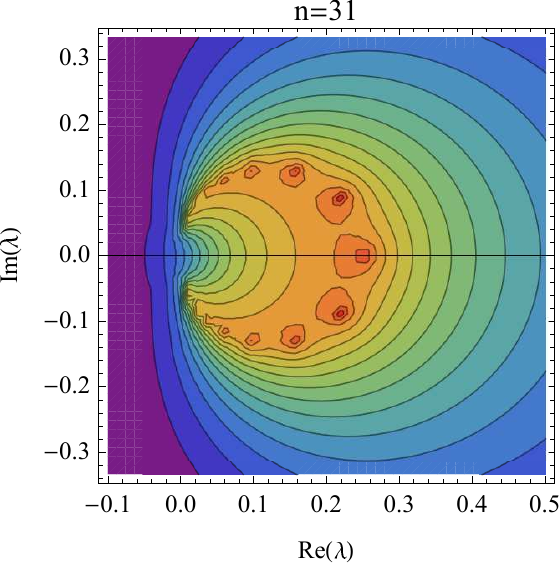}}\\
\subfloat{\includegraphics[width=0.45\linewidth]{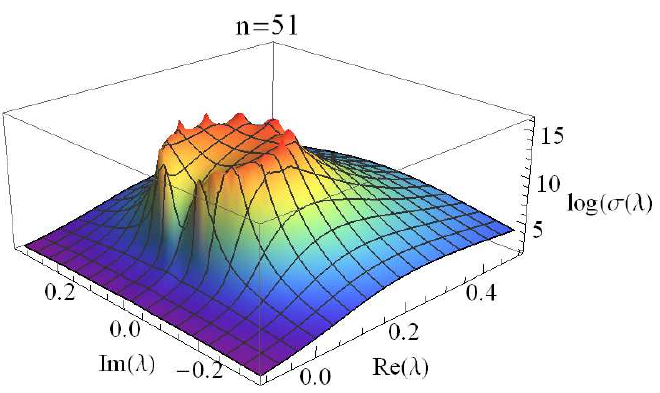}}
\subfloat{\includegraphics[width=0.4\linewidth]{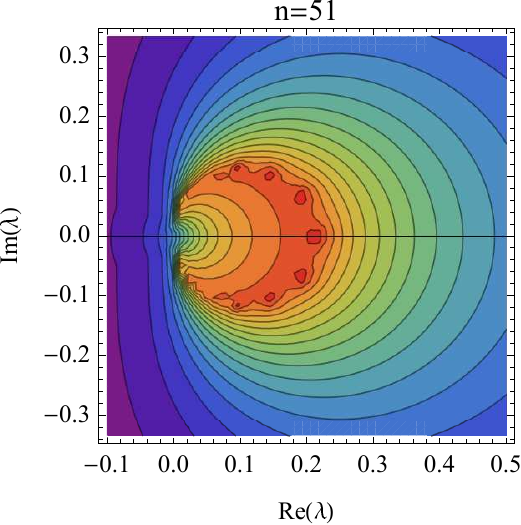}}\\
\subfloat{\includegraphics[width=0.45\linewidth]{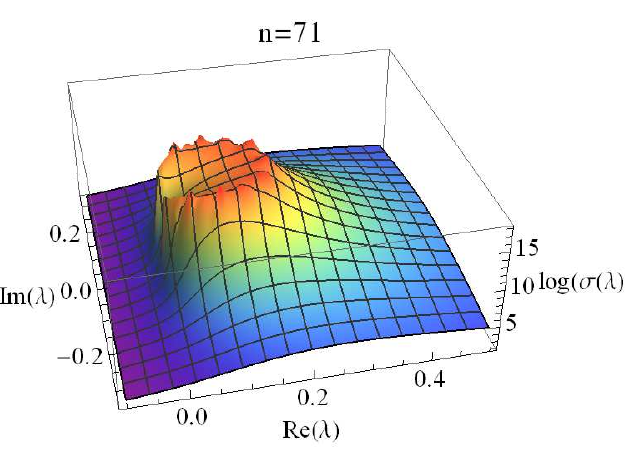}}
\subfloat{\includegraphics[width=0.4\linewidth]{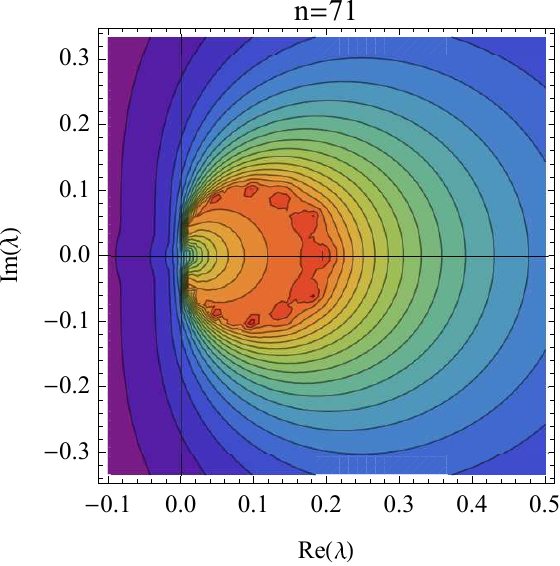}}
\caption{Resolvent of \(B^+\) for Chebyshev polynomials \(T_n(x)\).}
\label{Fig:15-20}
\end{figure}

\begin{expt} {\textbf{ CPC for Jacobi Polynomials}}

The Jacobi polynomials \(P_n^{\alpha ,\beta }(x)\), with \(-1\leq x\leq 1\) have the weight function \(w(x)=\xi (x)=(1-x)^{\alpha }(1+x)^{\beta }\)
with

\begin{equation*}
M_k=\int _{-1}^1(1-x)^{\alpha }(1+x)^{\beta } x^k dx<\infty , \forall k\geq 1 \text{and } \alpha ,\beta \geq 0.
\end{equation*}

Numerical results for the moments \(M_k\) are shown in Fig. 10.

\begin{figure}[H]
\centering
\includegraphics[scale=1.0]{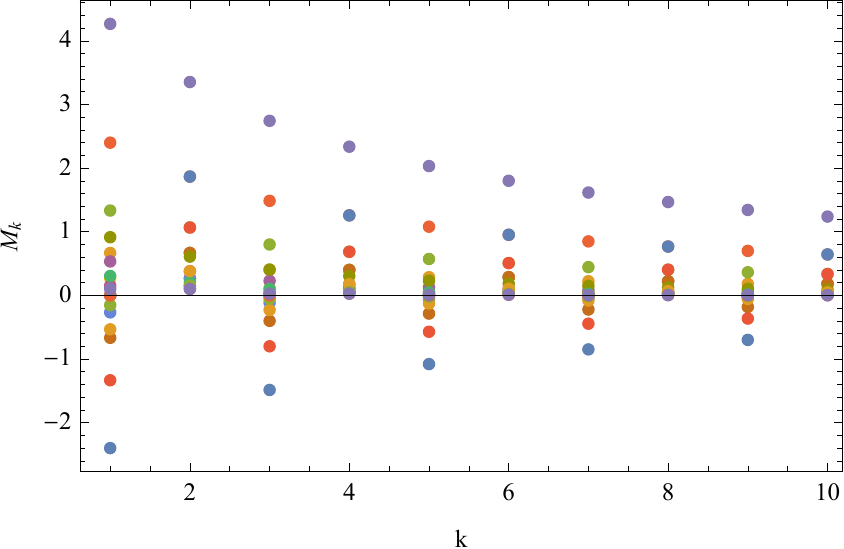}
\caption{The condition \(M_k\) for \(k=\)1,2,...,10 using \(\xi (x)=(1-x)^{\alpha }(1+x)^{\beta }\) for \(\alpha =0,1,2,3\) and \(\beta =1,2,3,4\).}
\label{Fig:21}
\end{figure}

In Fig. 11, the eigenvalues of the matrices \(B^+\) are presented. We used Jacobi polynomial \(P_n^{2,2}(x)\) with \(n=10:10:100\). For \(\alpha
=\beta -2\), we used \(\xi (x)=(1-x)^2(1+x)^2\) which is a symmetric function on \([-1,1]\). In Fig. 12, the resolvent of the matrix \(B^+\) are
presented.

\begin{figure}[H]
\centering
\includegraphics[scale=1.0]{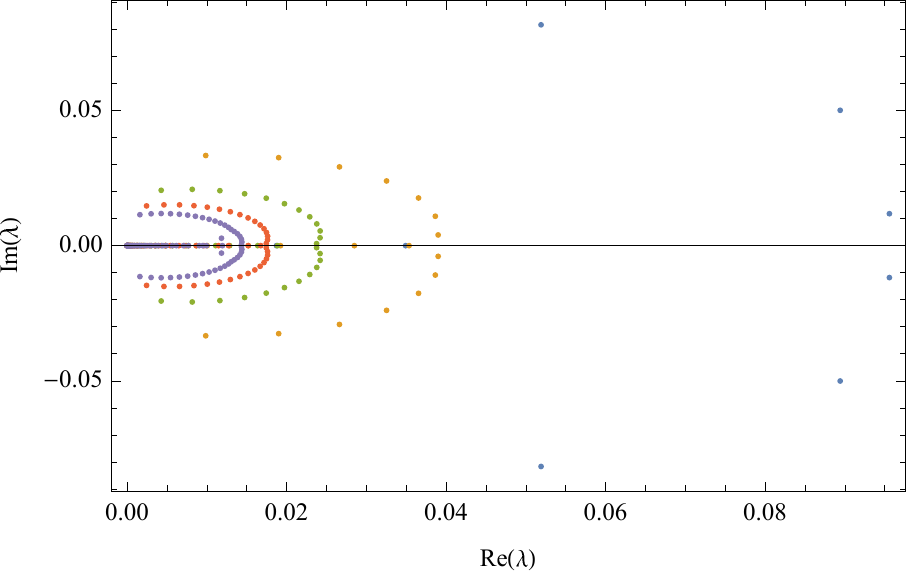}
\caption{Eigenvalues of \(B^+\) using Jacobi polynomial \(P_n^{2,2}(x)\), with \(n=10,30,50,70,90\).}
\label{Fig:22}
\end{figure}

\begin{figure}
\subfloat{\includegraphics[width=0.45\linewidth]{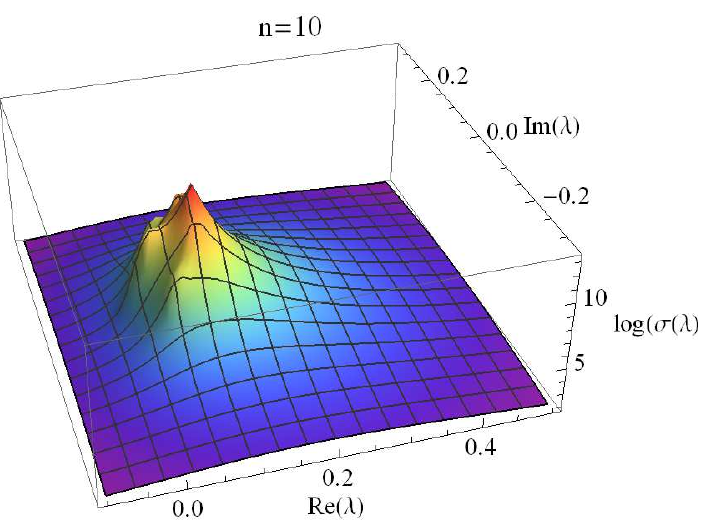}} 
\subfloat{\includegraphics[width=0.4\linewidth]{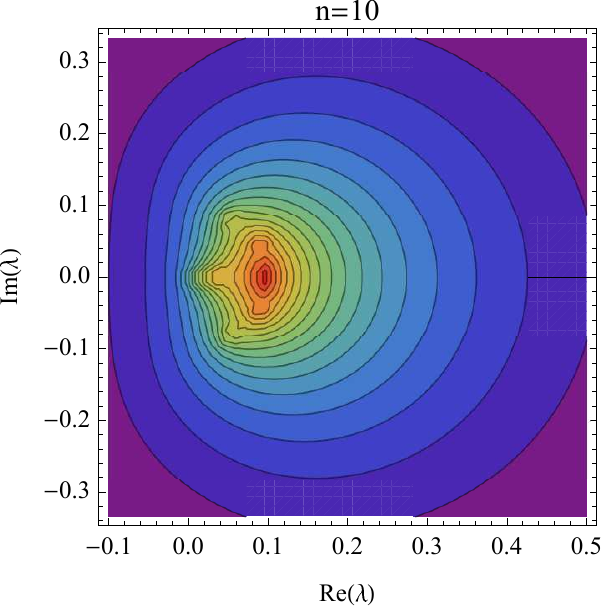}}\\
\subfloat{\includegraphics[width=0.45\linewidth]{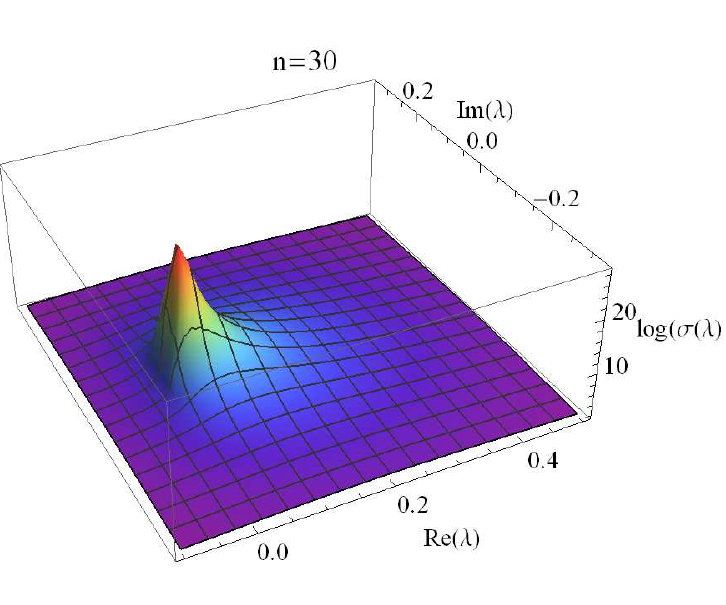}}
\subfloat{\includegraphics[width=0.4\linewidth]{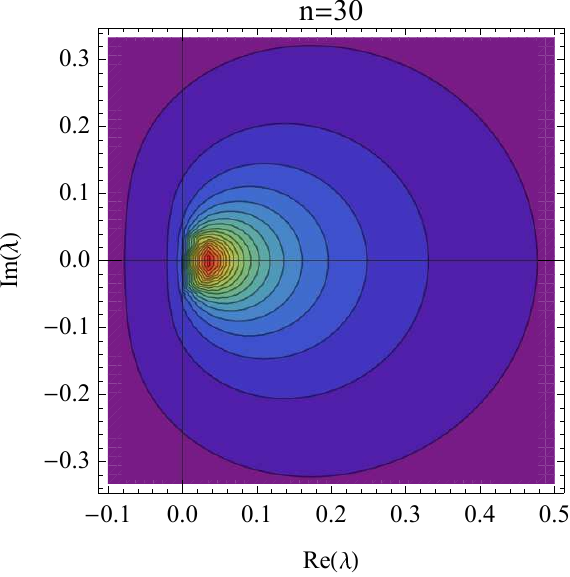}}\\
\subfloat{\includegraphics[width=0.45\linewidth]{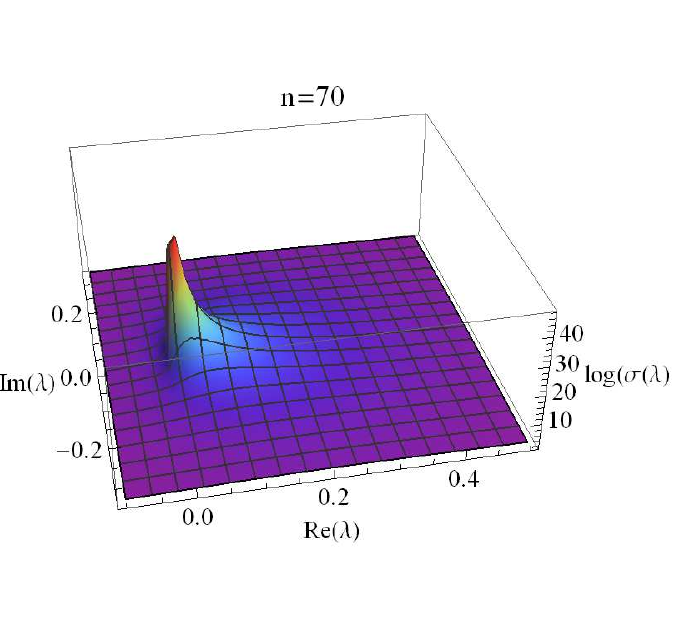}}
\subfloat{\includegraphics[width=0.4\linewidth]{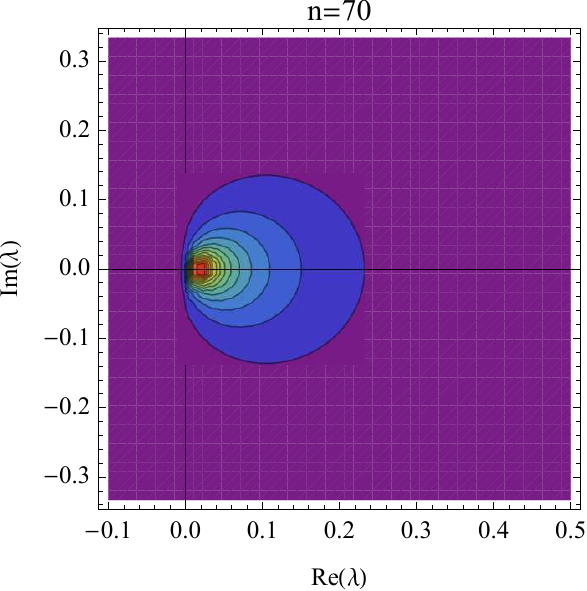}}
\caption{Resolvent of \(B^+\) using Jacobi polynomial\(P_n^{2,2}(x)\).}
\label{Fig:23-28}
\end{figure}

In Fig. 13, the eigenvalues of the matrices \(B^+\) based on Jacobi polynomials are presented. We used Jacobi polynomial \(P_n^{1,0}(x)\) with \(n=10,30,\text{... },90\). In this case we use \(\xi (x)=w(x)=1-x\). , which is not symmetric on \([-1,1]\). The figure and the tests show that the CPC is always satisfied. For the structure of the resolvent in this case, see Fig. 14.
\end{expt}

\begin{figure}[H]
\subfloat{\includegraphics[width=0.55\linewidth]{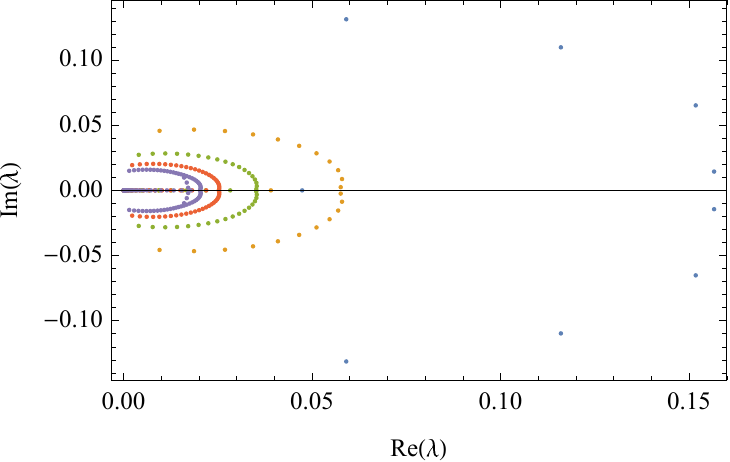}} 
\subfloat{\includegraphics[width=0.55\linewidth]{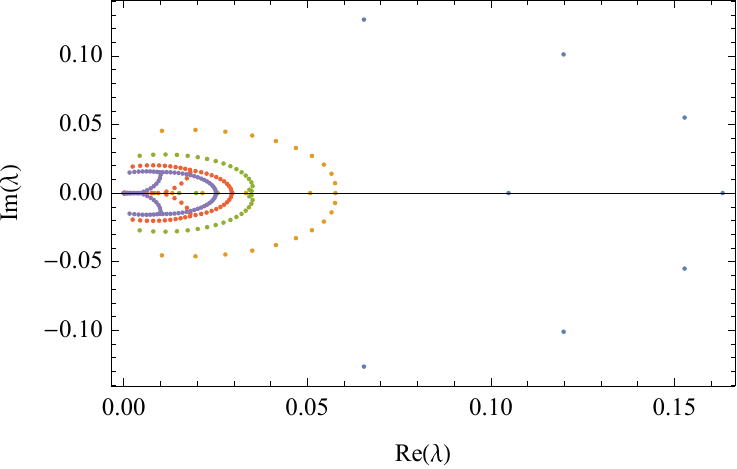}}
\caption{Eigenvalues of \(B^+\),left panel, and \(B^-\), right panel, using Jacobi polynomial \(P_n^{1,0}(x)\),
with \(n=10,30,50,70,90\).}
\label{Fig:29-30}
\end{figure}

\begin{figure}
\subfloat{\includegraphics[width=0.45\linewidth]{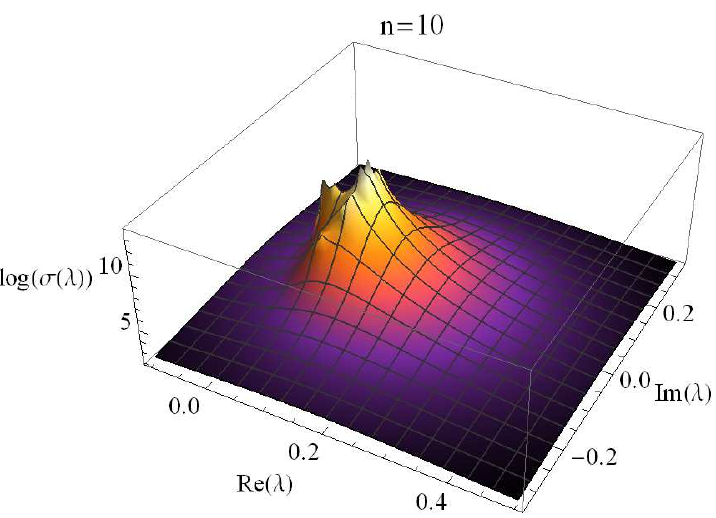}} 
\subfloat{\includegraphics[width=0.4\linewidth]{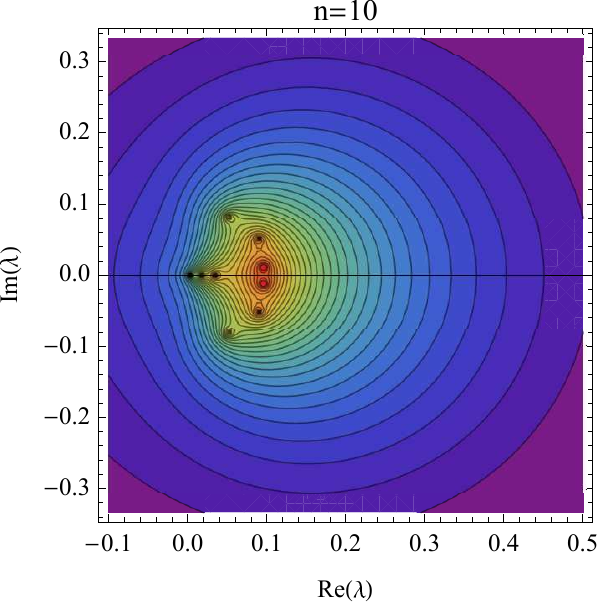}}\\
\subfloat{\includegraphics[width=0.45\linewidth]{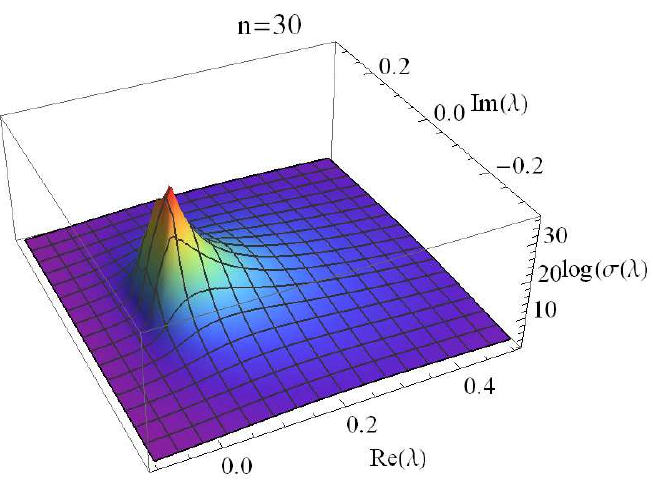}}
\subfloat{\includegraphics[width=0.4\linewidth]{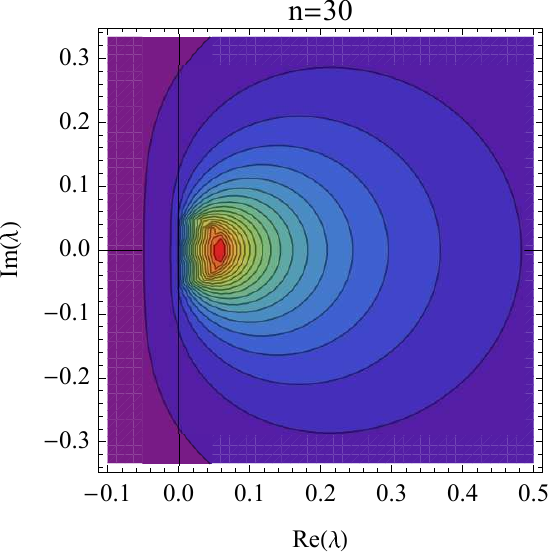}}\\
\subfloat{\includegraphics[width=0.45\linewidth]{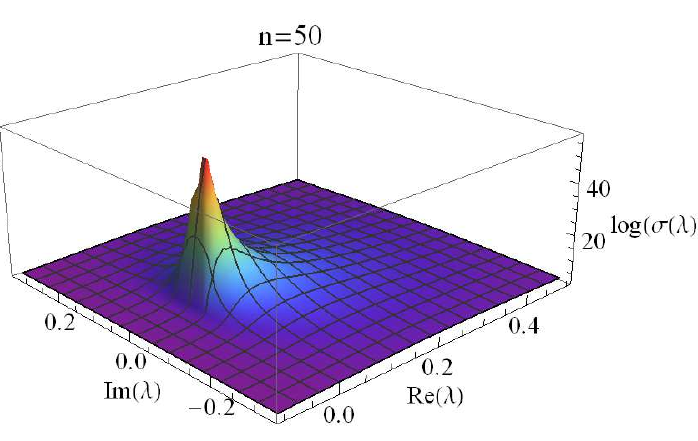}}
\subfloat{\includegraphics[width=0.4\linewidth]{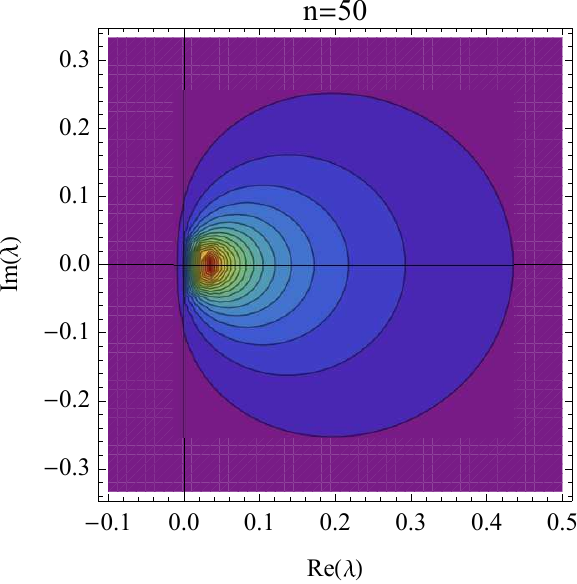}}
\caption{Resolvent of \(B^+\) using Jacobi polynomial\(P_n^{1,0}(x)\).}
\label{Fig:31-36}
\end{figure}

					
\begin{expt}{\textbf{ CPC for Gegenbauer Polynomials}}

Another generalization of both Legendre and Chebyshev polynomials are Gegenbauer polynomials \(C_n^{\eta }(x)\). { }The weight function for Gegenbauer
polynomials is \(\xi (x)=w(x)=\left(1-x^2\right)^{\eta -1/2}\) which is symmetric on \([-1,1]\) and satisfying

\begin{equation*}
M_k=\int _{-1}^1\left(1-x^2\right)^{\eta -1/2} x^kdx\geq 0, \forall k\geq 1 \text{and } \eta \geq 0.
\end{equation*}

The moments \(M_k\) are shown in Fig. 15 for different values \(\eta\).

\begin{figure}[H]
\centering
\includegraphics[scale=1.0]{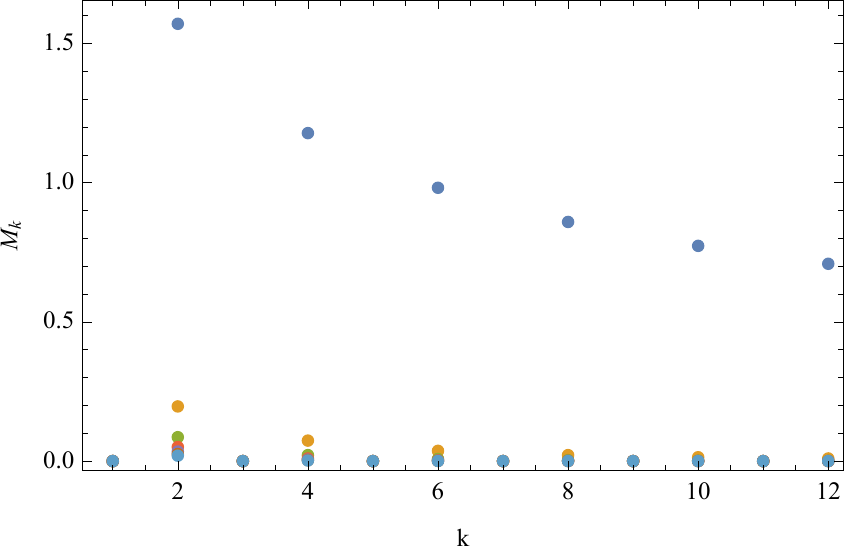}
\caption{The moments \(M_k\) for \(k=\)1,2,...,10 using \(\xi (x)=\left(1-x^2\right)^{\eta -\frac{1}{2}}\) for \(\eta =0,2,4,\ldots,12\).}
\label{Fig:37}
\end{figure}

Figure 16 shows eigenvalues for two different \(\eta\), \(\eta =2\) and \(\eta =10\) which satisfy the CPC. In Figs. 17 and 18 the corresponding resolvent structures are shown.

\begin{figure}[H]
\subfloat{\includegraphics[width=0.55\linewidth]{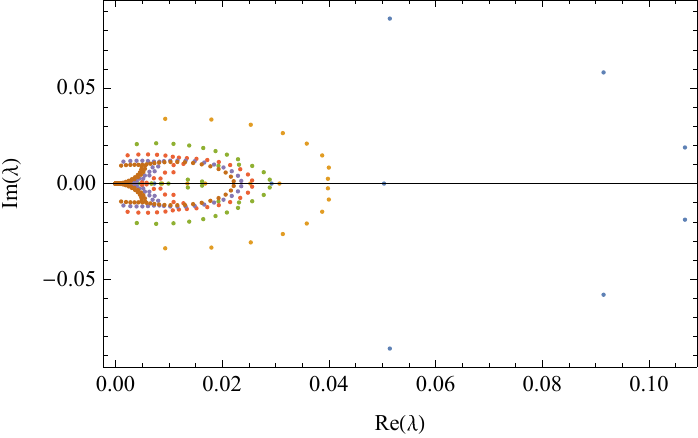}} 
\subfloat{\includegraphics[width=0.55\linewidth]{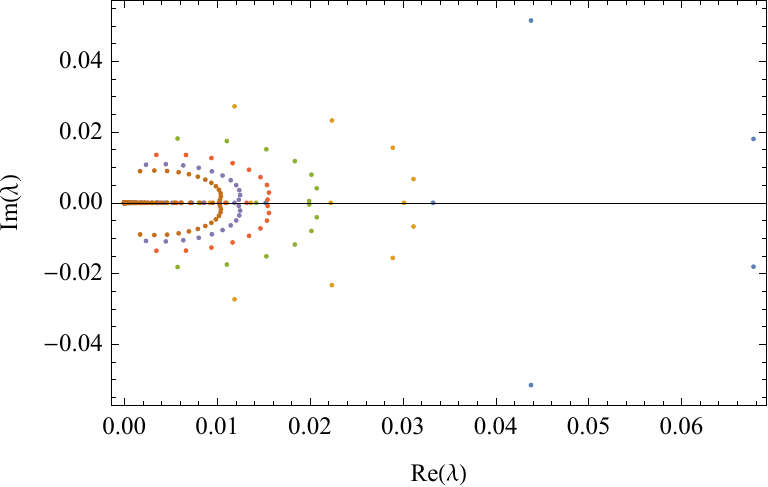}}
\caption{Eigenvalues of \(B^+\) using Gegenbauer polynomial \(C_n^2(x)\), left panel, and \(C_n^{10}(x)\), right panel. For both cases \(n=\{10,20,\ldots,100\}\).}
\label{Fig:38-39}
\end{figure}

\begin{figure}[H]
\subfloat{\includegraphics[width=0.5\linewidth]{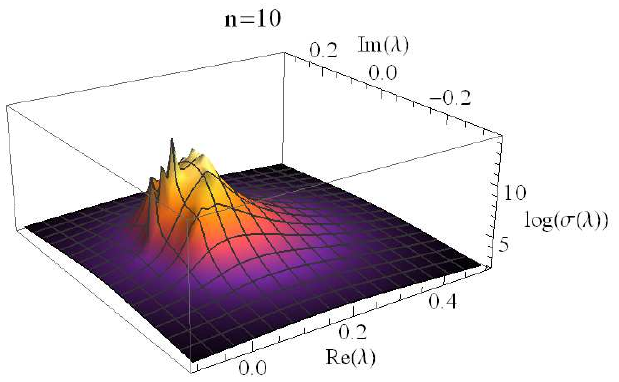}} 
\subfloat{\includegraphics[width=0.45\linewidth]{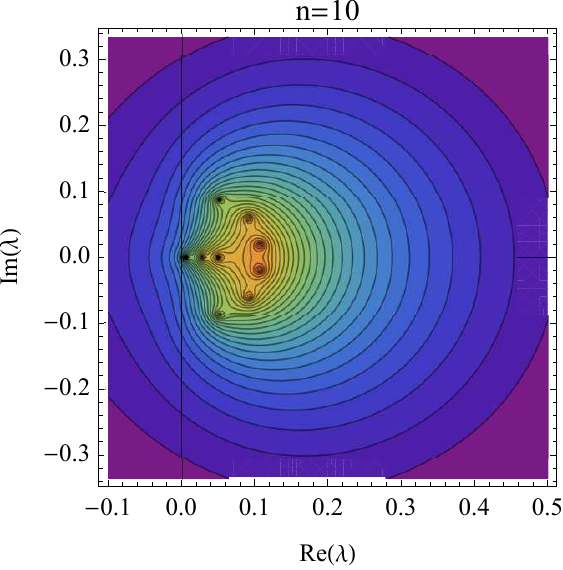}}
\caption{Resolvent of \(B^+\) using Gegenbauer polynomial \(C_n^2(x)\).}
\label{Fig:40-41}
\end{figure}

\begin{figure}[H]
\subfloat{\includegraphics[width=0.5\linewidth]{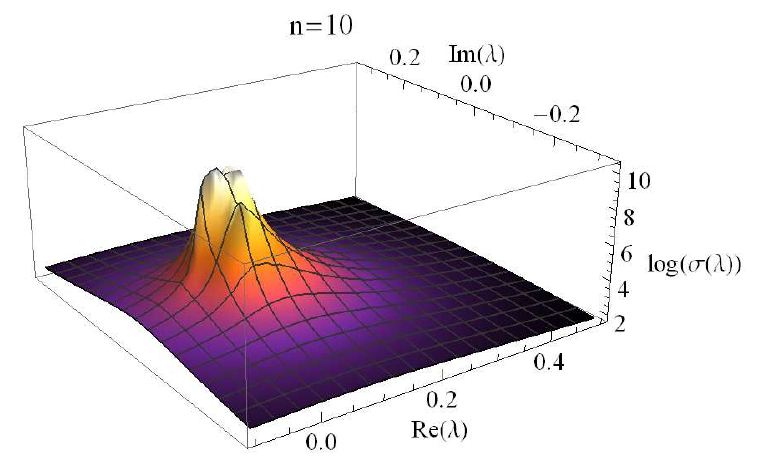}} 
\subfloat{\includegraphics[width=0.45\linewidth]{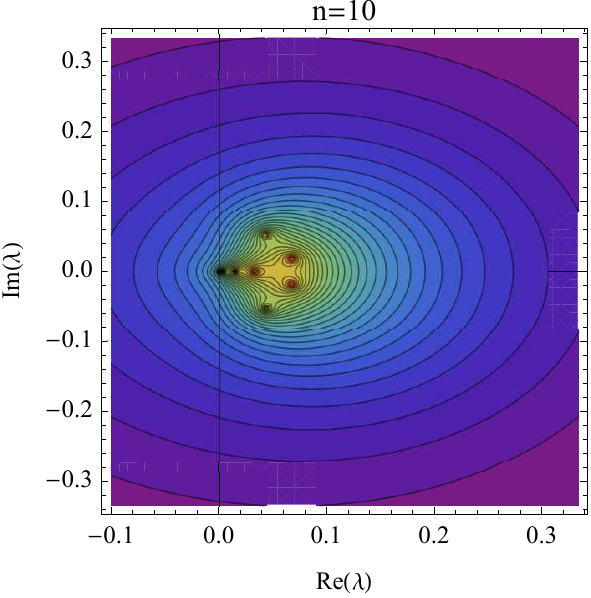}}
\caption{Resolvent of \(B^+\) using Gegenbauer polynomial \(C_n^{10}(x)\).}
\label{Fig:40-41}
\end{figure}

\end{expt}


\begin{expt}{\textbf{ CPC for Laguerre}}

In this experiment, we verify the CPC conjecture for Laguerre polynomials \(L_n(x)\) defined on \(\mathbb{R}^+=(0,\infty )\) with weight function
\(\xi (x)=w(x)=e^{-x}\). The moments \(M_k\) are shown in Fig. 19 as given by

\begin{equation*}
M_k=\int _0^{\infty }e^{-x} x^kdx\geq 1, \forall k\geq 1.
\end{equation*}

\begin{figure}[H]
\centering
\includegraphics[scale=1.0]{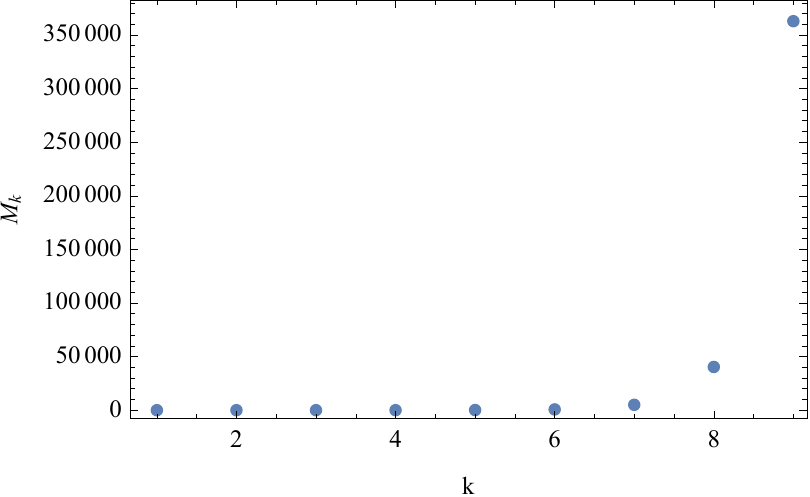}
\caption{The moments \(M_k\) for \(k=\)1,2,...,10 using \(\xi (x)=e^{-x}\).}
\label{Fig:44}
\end{figure}

Next, we verify the CPC for Laguerre polynomials \(L_n\) using \(n=10,20,\text{... },60\). The calculations are given in Fig. 20 for eigenvalues and in Fig. 21 for the resolvent.

\begin{figure}[H]
\subfloat{\includegraphics[width=0.55\linewidth]{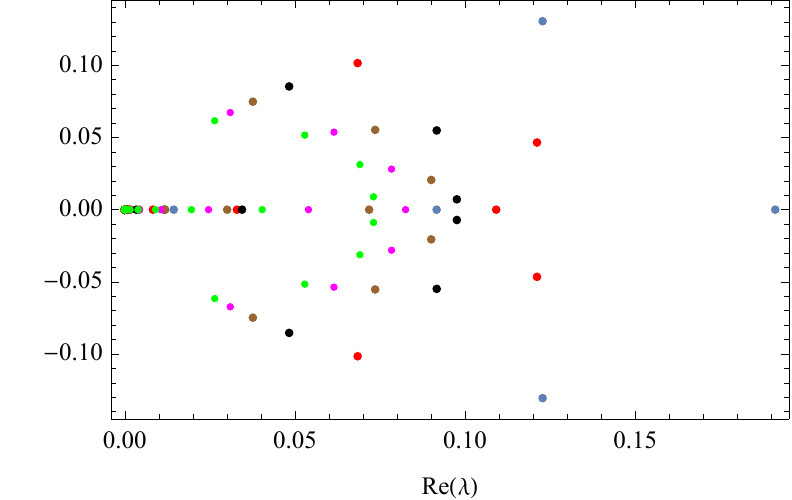}} 
\subfloat{\includegraphics[width=0.55\linewidth]{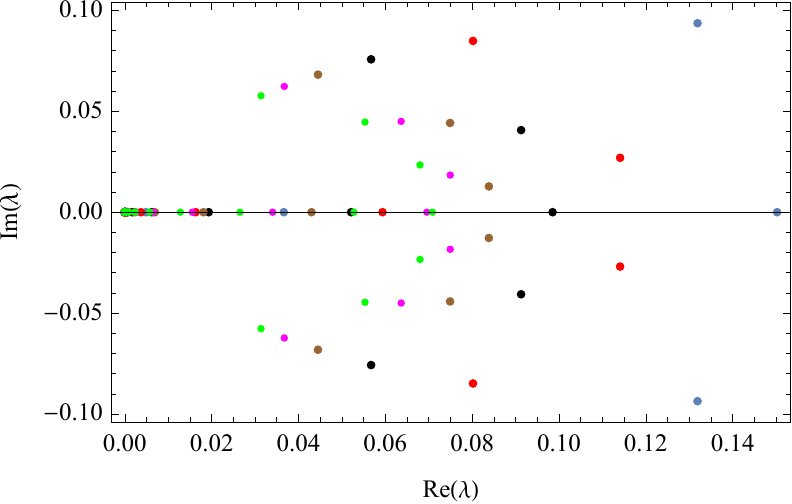}}
\caption{Eigenvalues of \(B^+\), right, and \(B^-\), left, using \(L_n(x)\) with \(n=10,20,\ldots,60\).}
\label{Fig:45-46}
\end{figure}

\begin{figure}
\subfloat{\includegraphics[width=0.45\linewidth]{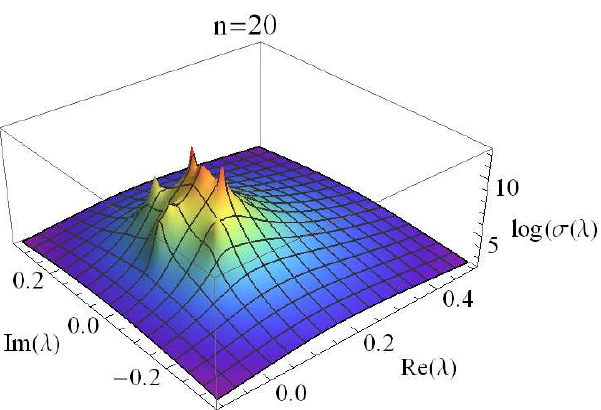}} 
\subfloat{\includegraphics[width=0.4\linewidth]{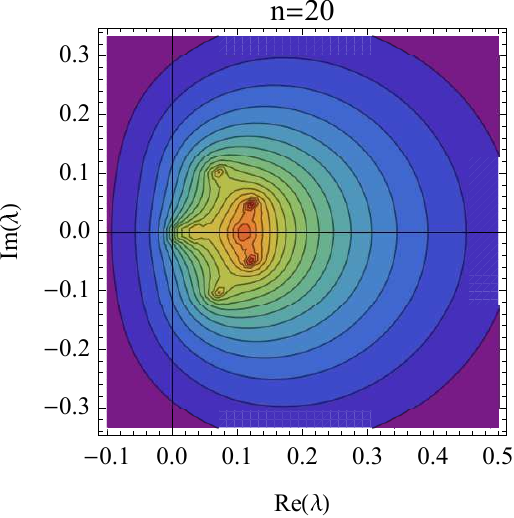}}\\
\subfloat{\includegraphics[width=0.45\linewidth]{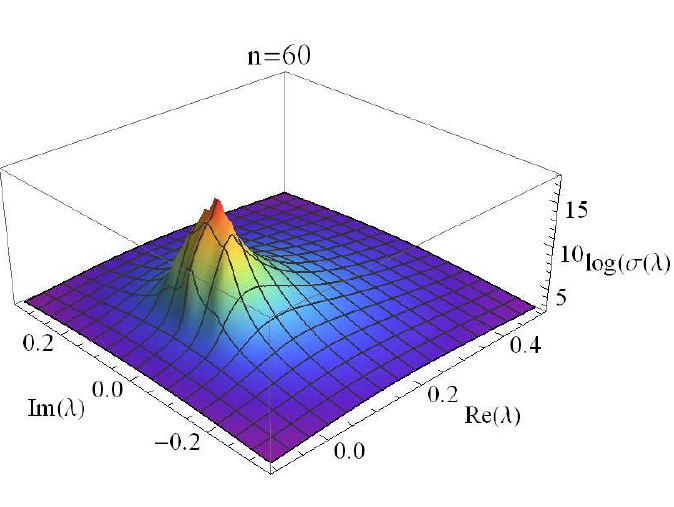}}
\subfloat{\includegraphics[width=0.4\linewidth]{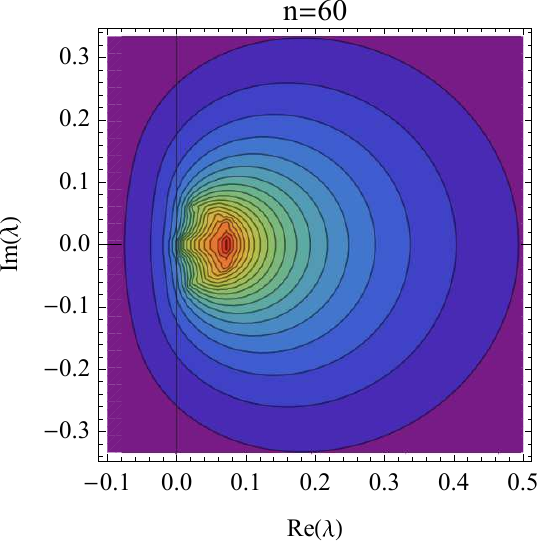}}
\caption{Resolvent of \(B^+\) using Laguerre polynomial \(L_n(x)\).}
\label{Fig:47-50}
\end{figure}

\end{expt}

\newpage

\begin{expt}{\textbf{ CPC for Hermit Polynomials}}

In this section, we verify the CPC conjecture for Hermite polynomials \(H_n(x)\) defined on \(\mathbb{R}=(-\infty ,\infty )\) with \(w(x)=e^{-x^2}\).
In this case the moments \(M_k\) are shown in Fig. 22 and computed by

\begin{equation*}
M_k=\int _{-\infty }^{\infty }e^{-x^2} x^kdx\geq 0, \forall k\geq 1.
\end{equation*}

\begin{figure}[H]
\centering
\includegraphics[scale=1.0]{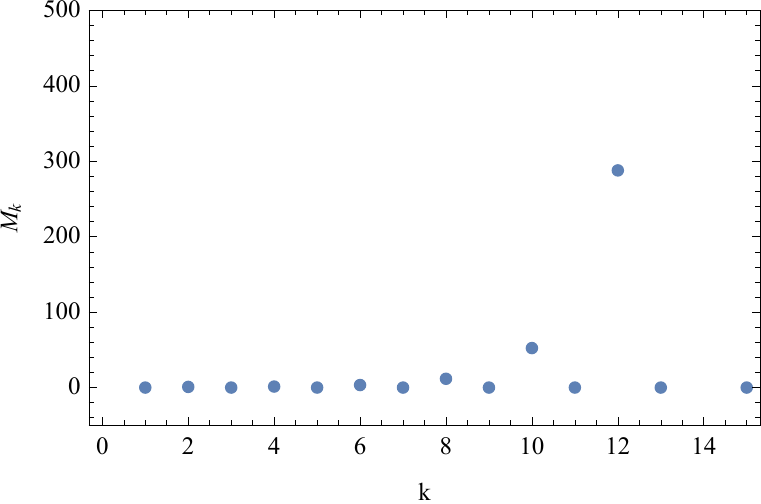}
\caption{The moments \(M_k\) for Hermit polynomials with \(k=1,2,\text{...},15.\) using \(\xi (x)=e^{-x^2}\).}
\label{Fig:51}
\end{figure}

The calculations of the spectrum of \(B^{\pm }\) are given in Fig. 23 and Fig. 24.

\begin{figure}[H]
\subfloat{\includegraphics[width=0.55\linewidth]{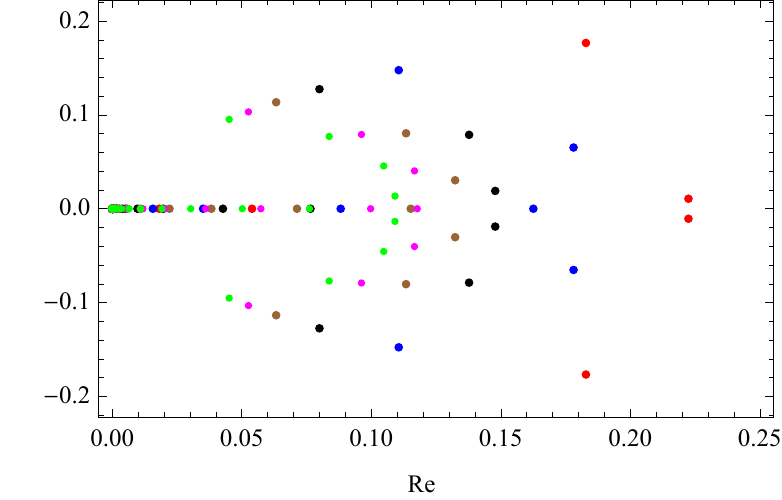}} 
\subfloat{\includegraphics[width=0.55\linewidth]{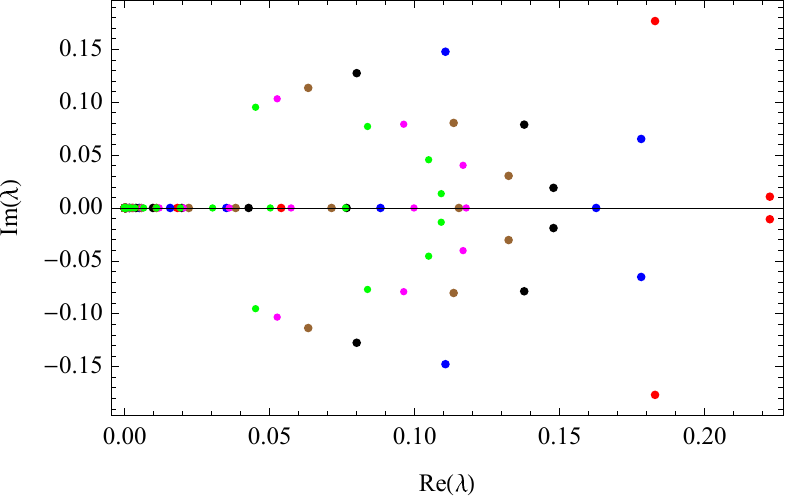}}
\caption{Eigenvalues of \(B^+\), right, and \(B^-\), left, using \(H_n(x)\) with \(n=10,20,\ldots,60\).}
\label{Fig:52-53}
\end{figure}

\begin{figure}
\subfloat{\includegraphics[width=0.45\linewidth]{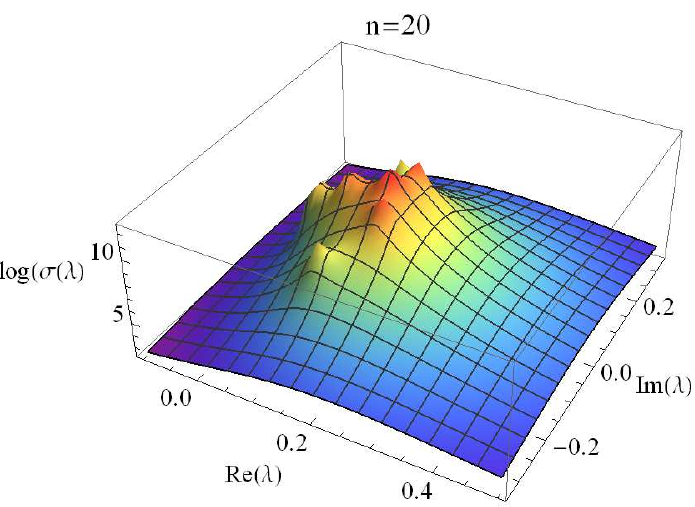}} 
\subfloat{\includegraphics[width=0.4\linewidth]{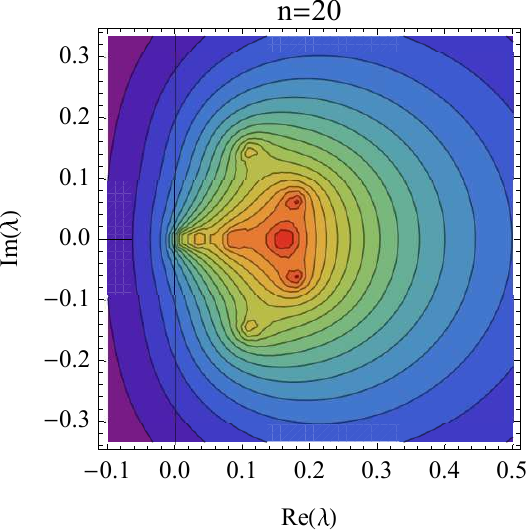}}\\
\subfloat{\includegraphics[width=0.45\linewidth]{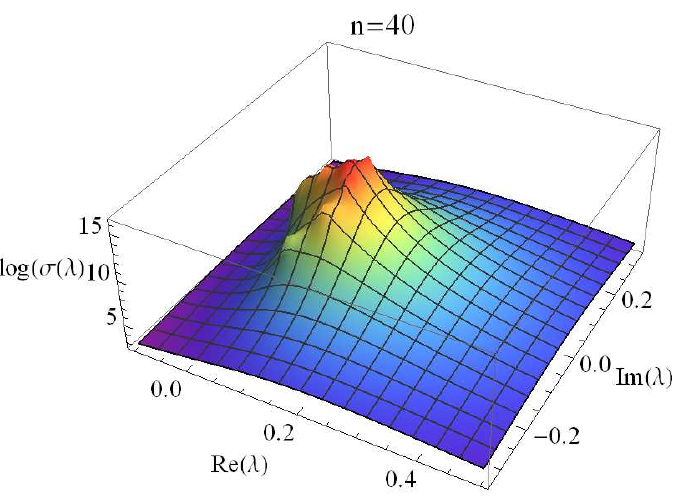}}
\subfloat{\includegraphics[width=0.4\linewidth]{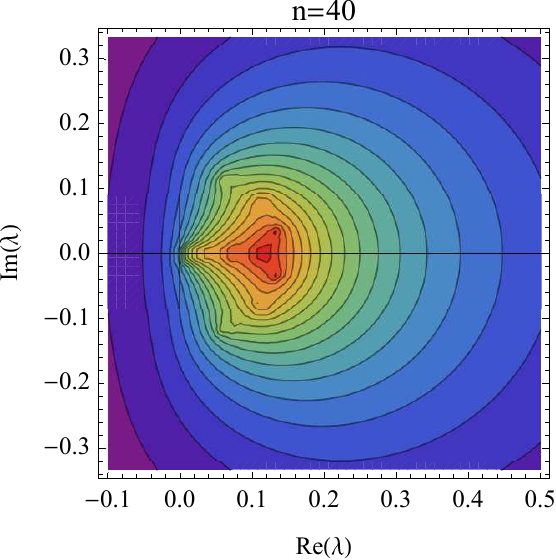}}
\caption{Resolvent of \(B^+\) using Hermit polynomials.}
\label{Fig:54-57}
\end{figure}

\end{expt}


\subsection{Verification of NPC}

In this section we discuss the case of \(\xi (x)\neq w(x)\). Many functions \(\xi (x)\) can be used, we are mainly interested in the case \(\xi (x)=1\),
which is corresponding to the approximation defined in (\ref{equation:1}) and (\ref{equation:5}). In this case the sufficient condition is not satisfied, which means no guarantee
for the correctness of the conjecture. We will see later that for some cases the conjecture will be satisfied while in other cases it will not.

\begin{expt}{\textbf{ NPC for Chebyshev polynomials}}

In this experiment we consider the matrices \(B^{\pm }\) defined in NPC conjecture using Chebyshev polynomials and \(\xi (x)=1\). The eigenvalue
calculations for \(B^+\) are shown in Fig. 25. From Fig. 25, it is clear that the real part of the eigenvalues are positive. 

\begin{figure}[H]
\subfloat{\includegraphics[width=0.55\linewidth]{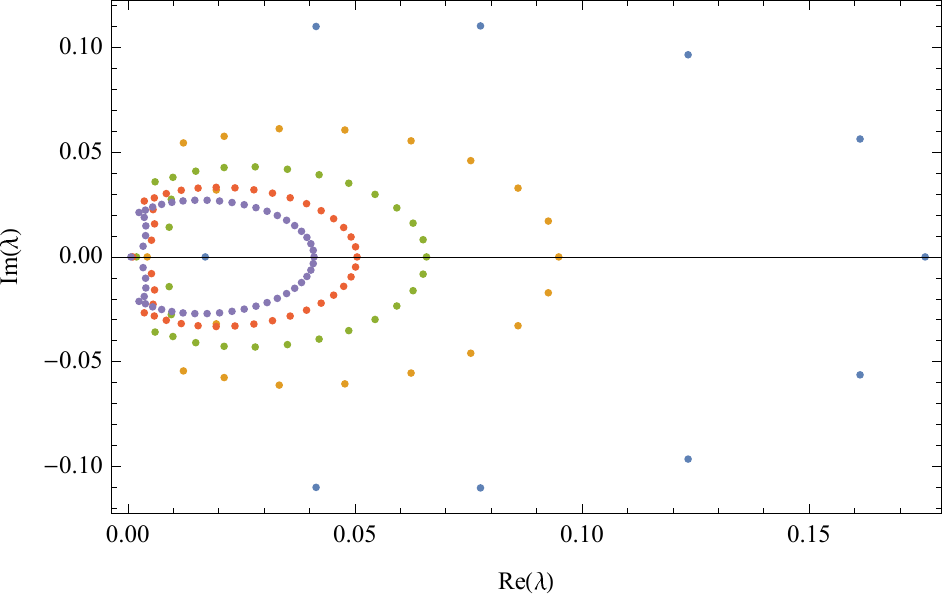}} 
\subfloat{\includegraphics[width=0.55\linewidth]{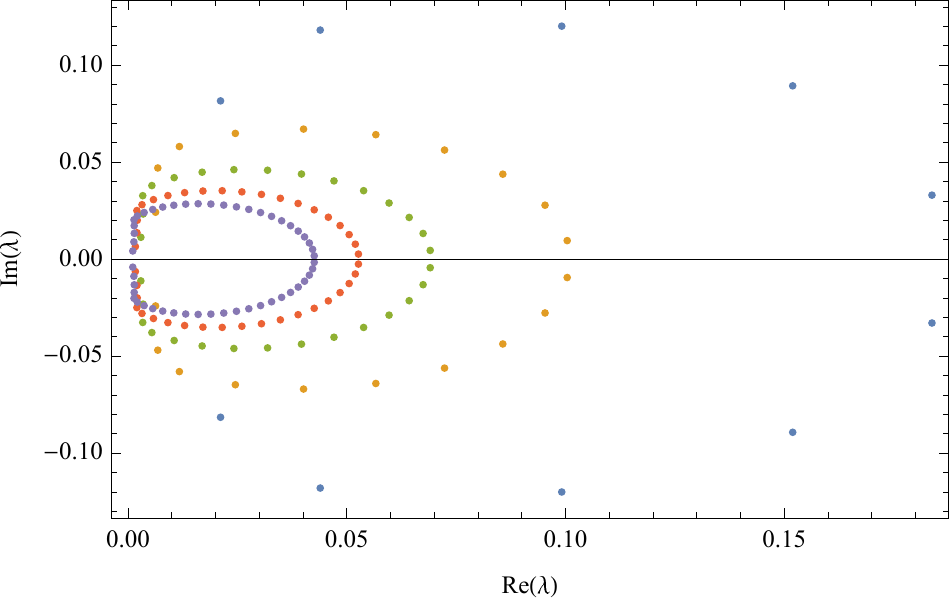}}
\caption{NPC for Chebyshev polynomials \(T_n(x)\), left, and \(U_n(x)\), right, \(n=10,20, \ldots ,50\).}
\label{Fig:58-59}
\end{figure}
\end{expt}


\begin{expt}{\textbf{ NPC for Jacobi polynomials}}

In this experiment we verify NPC for Jacobi Polynomials \(P_n^{\alpha ,\beta }(x)\). NPC is not valid for all \(\alpha\) and \(\beta\). For the \(\alpha
=\beta =2\), the NPC is not verified, see Fig. 26. For \(\alpha =1\) and \(\beta =0\), the conjecture is true, see Fig. 27.

\begin{figure}[H]
\centering
\includegraphics[scale=1.0]{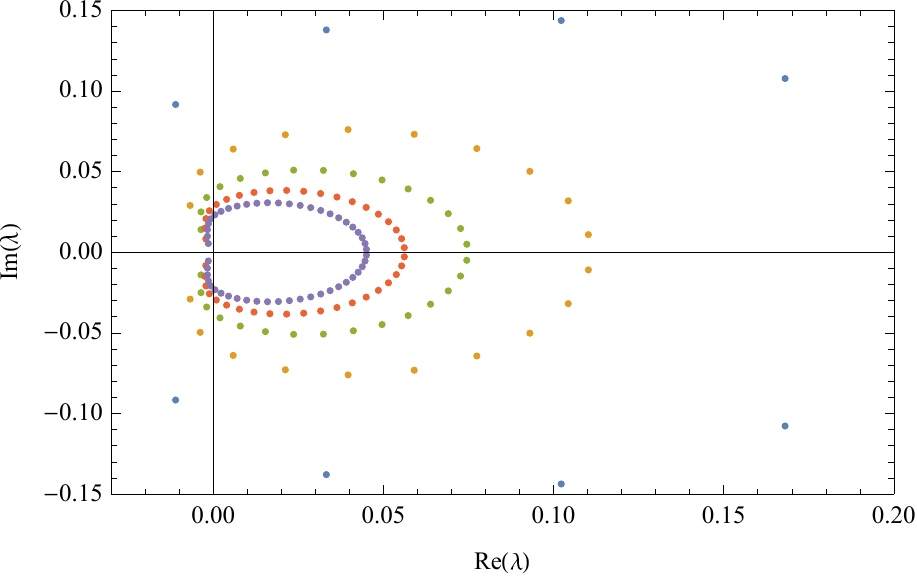}
\caption{NPC for Jacobi Polynomials \(P_n^{2,2}(x)\), \(n=10,20,\ldots,50\).}
\label{Fig:60}
\end{figure}

\begin{figure}[H]
\centering
\includegraphics[scale=1.0]{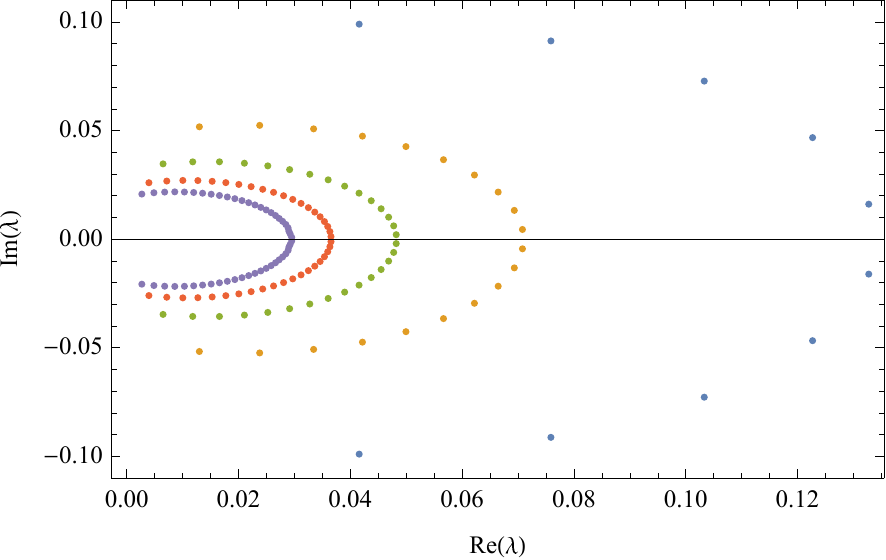}
\caption{NPC for Jacobi Polynomials \(P_n^{1,0}(x)\), \(n=10,20,\ldots,50\).}
\label{Fig:61}
\end{figure}

\end{expt}


\begin{expt}{\textbf{ NPC for Gegenbauer polynomials}}

For Gegenbauer polynomials with \(\xi (x)=1\) and \(\eta =2\) (Fig. 28) and \(\eta =10\) (Fig. 29) the NPC is not satisfied.

\begin{figure}[H]
\centering
\includegraphics[scale=1.0]{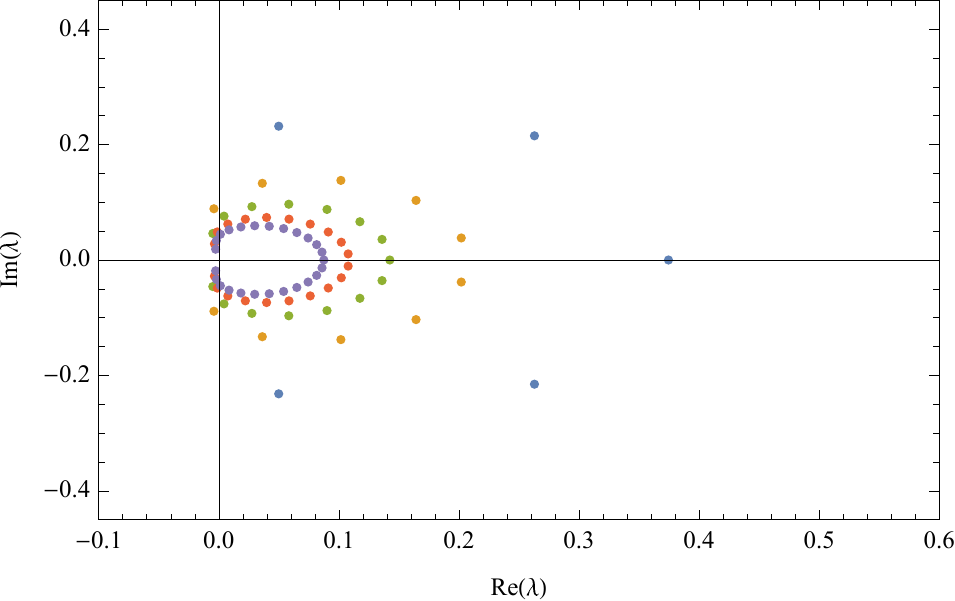}
\caption{NPC for Gegenbauer Polynomials \(C_n^2(x)\), \(n=5,10,15,20,25\).}
\label{Fig:62}
\end{figure}

\begin{figure}[H]
\centering
\includegraphics[scale=1.0]{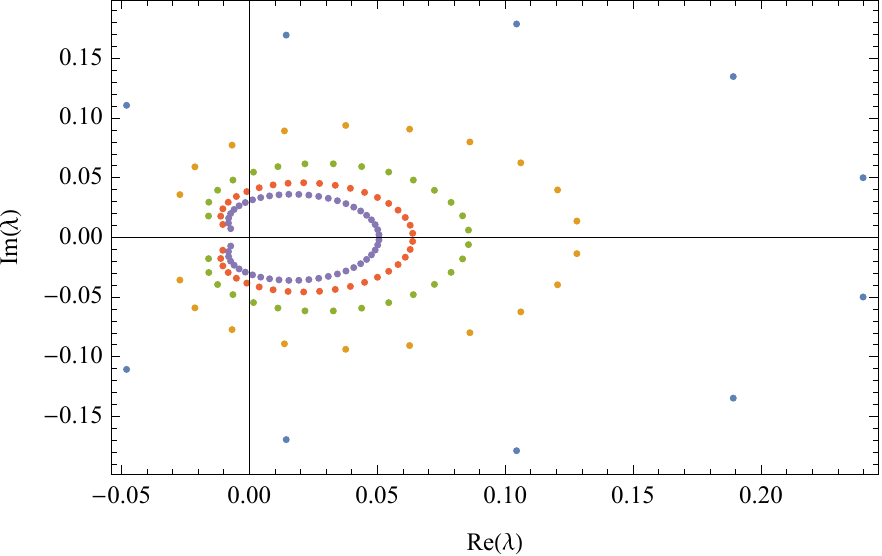}
\caption{NPC for Gegenbauer Polynomials \(C_n^{10}(x)\), \(n=10,20,\ldots,50\).}
\label{Fig:63}
\end{figure}

\end{expt}


\begin{expt}{\textbf{ NPC for Laguerre polynomials}}

The verifications of NPC for Laguerre polynomials are given in Fig. 30 and Fig. 31. Although the calculations show that \(B^+\) has positive eigenvalues
the matrix \(B^-\) contrary shows negative eigenvalues.
\begin{figure}[H]
\centering
\includegraphics[scale=1.0]{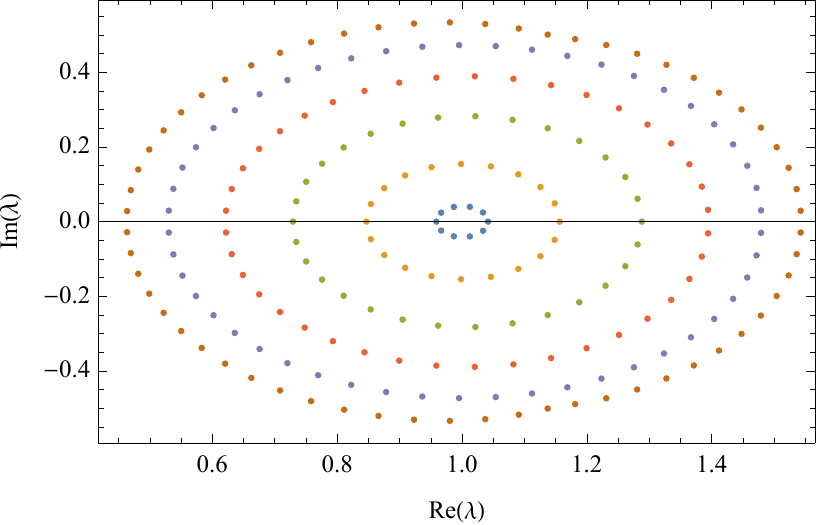}
\caption{Eigenvalues of \(B^+\) using \(L_n(x)\) with \(n=10,20, \ldots ,50\).}
\label{Fig:64}
\end{figure}

\begin{figure}[H]
\centering
\includegraphics[scale=1.0]{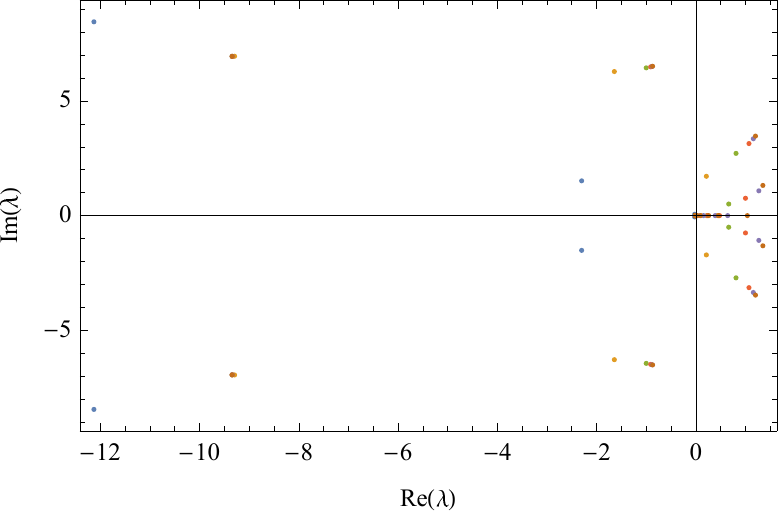}
\caption{Eigenvalues of \(B^-\) using \(L_n(x)\) with \(n=5,15,20,25\).}
\label{Fig:65}
\end{figure}

\end{expt}


\begin{expt}{\textbf{ NPC for Hermit polynomials}}

Both for matrices \(B^{\pm }\) we find negative real parts of the eigenvalues. We will skip the plots here.
\end{expt}


\subsection{Poly-Sinc Matrices}

Sinc points are related to a family of orthogonal functions, Sinc cardinal functions. Although these points are not roots of orthogonal polynomials,
but they are used in two effective approximations, Sinc approximation and Poly-Sinc approximation. For Sinc approximation, a similar matrix conjecture
has been formulated by Stenger and Proved by Han and Xi in 2014 \cite{Han_2014}. In this section, we verify the NPC for Poly-Sinc polynomial matrices. Poly-Sinc
approximation on finite intervals based on the use of Sinc points as interpolation points in Lagrange approximation [3, 5]. This kind of Polynomial
approximation in connection with the conformal maps from the finite interval to $\mathbb{R}$. For the polynomial approximation defined in (\ref{equation:2}), we
have the basis functions

\begin{equation}
b_k(x)=\prod _{l=-M, l\neq k}^N \frac{x-x_l}{x_k-x_l}=\frac{v_m(x)}{\left(x-x_k\right)v'_m\left(x_k\right)}, m=M+N+1,
 \label{equation:31}
\end{equation}

where \(x_k=\phi ^{-1}(k h)\) are the Sinc points. Poly-Sinc approximations shows an exponential decaying rate similar to Sinc approximations, with
smaller Lebesgue constant \cite{Maha_Phd}. Now define the functions

\begin{equation}
v_m(x)=\prod _{l=-M}^N x-x_l=, m=M+N+1.
 \label{equation:32}
\end{equation}

The family of polynomials \(\left\{b_k(x)\right\}\) are defined for any finite interval \((a,b)\) with conformal map \(\phi (x):(a,b)\longrightarrow
\mathbb{R}\) defined as

\begin{equation}
\phi (x)=\text{Log }\left(\frac{x-a}{b-x}\right)
 \label{equation:33}
\end{equation}

 and with a set of Sinc points \(x_k=\phi ^{-1}(k h)\) defined as

\begin{equation}
x_k=\frac{b e^{k h}+a}{1+e^{k h}},
 \label{equation:34}
\end{equation}

 For these types of polynomials we verify the NPC with \(\xi (x)=1\). As a study interval we use we choose \([a,b]=[-1,1]\). It is known that Poly-Sinc
shows high accuracy even with used number of Sinc points \cite{Youssef_2016, Youssef_2019}. So, we will test the eigenvalues for \(B^{\pm }\) for not so huge numbers, roughly
we test up to 41 Sinc points. The result of these calculations are given n Fig. 32.

\begin{figure}[H]
\subfloat{\includegraphics[width=0.55\linewidth]{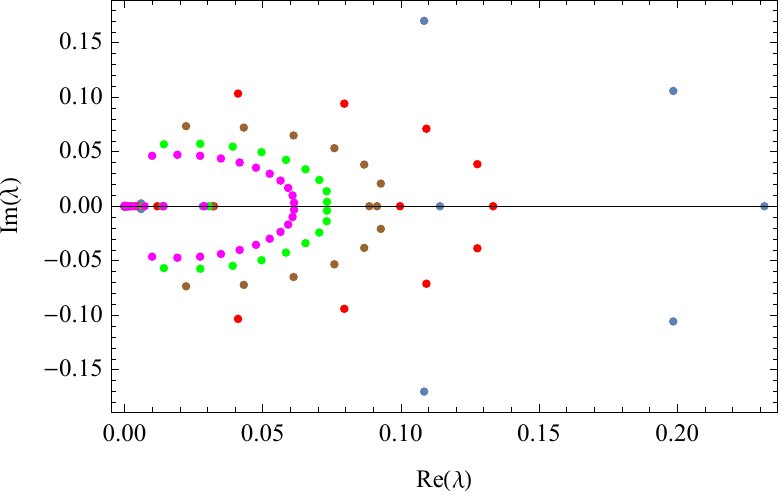}} 
\subfloat{\includegraphics[width=0.55\linewidth]{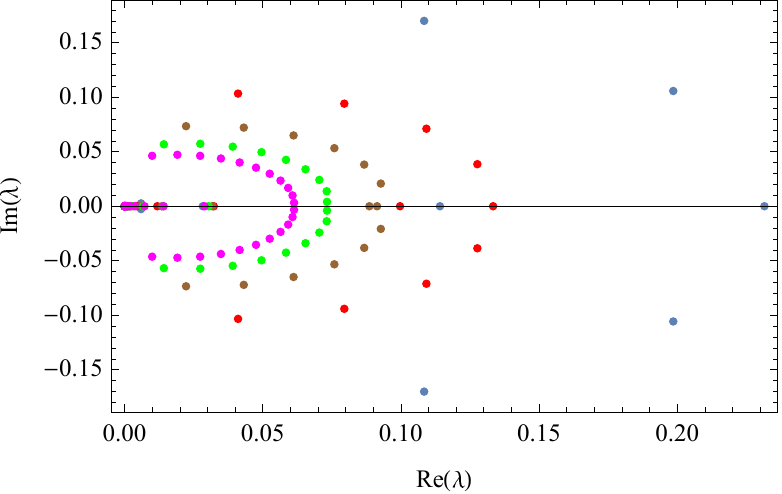}}
\caption{Eigenvalues of \(B^{\pm}\) using Sinc points on \([-1,1]\) with \(n=9,17,25,33,41\).}
\label{Fig66-67}
\end{figure}

If the finite interval \([a,b]\) is not \([-1,1]\), we define a one-to-one transformation from \([a,b]\) to \([-1,1]\). This transformation defines
a new distribution of Sinc points and maps the Lagrange basis to a new basis that satisfies the same properties as the old basis, for more details
see \cite{Stenger_2013}.


\section*{Conclusion}

In this paper we introduced a weighted form of polynomial Interpolation. The basis function are defined using sets of orthogonal polynomials and
their roots. As a result approximations of integral operators have been defined. The {``}New Polynomial Conjecture{''} has been verified/contracted
for set of orthogonal polynomials defined on finite, semi-infinite and infinite intervals. We introduce a reformulation of the conjecture to be verified
for all orthogonal polynomials. The numerical approach shows that this new conditioned conjecture is always true. Finally, we verified the conjecture
for a different set of polynomials called Poly-Sinc. Both Lagrange Interpolation using Sinc points or roots of orthogonal polynomials as interpolation
points yield exceptional rates of error for approximating the function and integral operators that are essential for the solution of PDEs.

\section*{Acknowledgement}

We are indebted to Frank Stenger for the discussions on these topics and his comments on the weighted approximation. The first author would like
to thank Bernd Kugelmann, University of Greifswald, for the fruitful discussions and comments during the preparation of the paper.


\begin{thebibliography}{14}

\bibitem{Stenger_86} 
F. Stenger, \textit{Explicit, Nearly Optimal, Linear Rational Approximations with Preassigned Poles}, Math. Comp. v.47, pp. 225-252, (1986) .
 
\bibitem{Stenger_2015} 
F. Stenger, G., Baumann, V.G., Koures, \textit{Computational Methods for Chemistry and Physics, and Schr{\" o}dinger in 3+1}. In: Sabin, J.R., Cabrera-Trujillo,
R. (eds.) Advances in quantum chemistry. Concepts of mathematical physics in chemistry: a tribute to Frank E. Harris, v. 71, pp. 265-298. Academic
Press, Amsterdam, (2015).
 
\bibitem{Stenger_2013} 
F. Stenger, M. Youssef, J. Niebsch, \textit{Improved approximation via use of transformations}, in: X. Shen, A.I. Zayed (Eds.), Multiscale Signal Analysis
and Modeling, Springer, NewYork, pp. 25-49, (2013).

\bibitem{Stenger_2011} 
F. Stenger, \textit{Handbook of Sinc Numerical Methods}, CRC Press, (2011).

\bibitem{Maha_Phd} 
M. Youssef, \textit{Poly-Sinc Approximation Methods}, PhD thesis, Math. Dept. German University in Cairo, (2017).

\bibitem{Meinardus_67} 
G. Meinardus, \textit{Approximation von Funktionen und ihre numerische Behandlung}, Springer Verlag, (1964).

\bibitem{Bogaert_2012} 
I. Bogaert, B. Michiels, and J. Fostier, \textit{{``}O(1) Computation of Legendre polynomials and Gauss-Legendre nodes and weights for parallel computing}, Siam Journal on Scientific Computing, vol. 34, no. 3, pp. C83-C101, (2012).

\bibitem{Gautschi_2009} 
W. Gautschi, \textit{How Sharp is Bernstein{'}s Inequality for Jacobi Polynomials?}, Department of Computer Sciences, Purdue University, (2009).
 
\bibitem{Zhao_2013} 
X. Zhao, L. Wang, Z. Xie, \textit{Sharp Error Bounds for Jacobi Expansions and Gegenbauer--Gauss Quadrature of Analytic Functions}, SIAM J. Numer. Anal.,
51(3), pp. 1443$--$1469, (2013).

\bibitem{Koornwinder_2018} 
T. Koornwinder, A. Kostenko, G. Teschl, \textit{Jacobi polynomials, Bernstein-type inequalities and dispersion estimates for the discrete Laguerre operator},
Advances in Mathematics, 333, pp. 796-821, (2018). 
 
\bibitem{Gautschi_2018} 
W. Gautschi and E. Hairer, \textit{On conjectures of Stenger in the theory of orthogonal polynomials}, E. J Inequal Appl, 2019: 159. https://doi.org/10.1186/s13660-019-2107-6, (2019).
 
\bibitem{Han_2014} 
L. Han, J. Xu, \textit{Proof of Stenger{'}s conjecture on matrix \(I^(-1)\) of Sinc methods}, Journal of Computational and Applied Mathematics 255, pp. 805-811,
(2014).

\bibitem{Youssef_2016} 
M. Youssef, G. Baumann, \textit{Collocation Method to Solve Elliptic Equations, Bivariate Poly-Sinc Approximation}, Journal of Progressive Research in Mathematics
(JPRM), ISSN: 2395-0218, 7(3), pp. 1079-1091 (2016).

\bibitem{Youssef_2019} 
M. Youssef, R. Pulch, \textit{Poly-Sinc Solution of Stochastic Elliptic Differential Equations}, http://arxiv.org/abs/1904.02017, (2019).


\end{thebibliography}
\end{document}